# Numerical approximation of systems modeling  large structure formation in cosmology

M.  Colombeau, Département de Mathématiques, Faculté des Sciences, Université des Antilles et de la Guyane, 97157, Pointe à Pitre Cedex, Guadeloupe (France).

## Abstract

Numerical approximations of two classical fluid dynamics systems modelling large structure formation in cosmology are proposed. These systems model nonrelativistic and relativistic fluids submitted to self-gravitation in an expanding background. They are obtained by an adaptation of an extension of the Godunov method, using delta wave projections, which was first introduced in Le Roux & al. for the system of pressureless fluid dynamics and which adapts to various systems of fluid dynamics.

## Introduction

In [1] the authors introduced the idea to project delta waves to obtain a Godunov type numerical scheme in pressureless fluid dynamics. The  method in [1] has been modified in Appendix 2 to treat the case of arbitrary changes in sign of velocity as needed in physical situations such as those considered in this paper. First the modified method is used for a numerical solution of the system modelling a Newtonian self-gravitating fluid in an expanding background with random initial conditions around a constant value for the energy density and around the value 0 for velocity. Thus velocity changes sign at random in the initial condition, and shows important changes in sign throughout the calculations. Then a similar scheme is presented in the relativistic case, as well as various multifluid extensions.

Either directly, or indirectly through the numerical technique of splitting of equations, this method applies to provide numerical schemes in numerous models in fluid dynamics [1,Appendix 2]. It is the purpose of this paper to put this method at the disposal of physicists, on the occasion of these two standard systems of equations of physics. The method is presented in 2D for convenience since the 2D calculations work on any PC. Its adaptation to 3D is straightforward. In order to obtain neat results after a few iterations one has choosen possibly non realistic values for the physical parameters (gravitation constant, velocity of light, expansion rate, values of the random irregularities of the initial medium).

In this paper we first consider the case of a single fluid modelled in the most usual way: continuity equation, Euler equation and Poisson equation, or their relativistic counterpart. The needed extension to multifluid flow in cosmology is straightforward since the cosmic fluids evolve according to their own pressure and are usually coupled only through gravitation ([10] p266, 268). Two multifluid extensions are given: relativistic+Newtonian (baryonic matter + dark matter between equivalence and decoupling), Newtonian+Newtonian (baryonic matter+dark matter after decoupling). One observes from 2D numerical simulations, that last 3 to 15 minutes on a standard PC, how the Newtonian and the relativistic systems can create structures in adequate conditions, the role of Jeans' length, the Meszaros effect, and the role of potential wells of dark matter at decoupling.

The expansion of the background is described by the "scale factor" a(t) which is a regular >0 increasing function of the time : a physical distance = unity at time 0 becomes a(t) at time t. The scale factor a(t) is obtained from the Friedman equations. One sets H(t)=$\dot{a}(t)/a(t)$ (= Hubble function). The function $\rho=\rho$(x,t) denotes the energy density and the function $\vec{u}=\vec{u}$(x,t) denotes the velocity vector; its components are stated (u,v,w), G = gravitation constant and $\Phi$= gravitation potential. Then the fluid is described by the following equations in comoving coordinates (i.e. spatial coordinates whose unit of length follows the expansion of the background: the spatial physical coordinates are obtained by multiplying the comoving coordinates by the scale factor a(t)). The time considered here is the usual physical time.

## 1. The equations in the Newtonian domain

They are the continuity equation, the Euler equation and the Poisson equation respectively (see [10] p233, [5] p294). The equations are usually stated with the state law p=0 (dark matter, baryonic matter after decoupling); one also considers the state law $p=\frac{1}{3}<v^2>\rho$, where $<v^2>$ is the mean squared velocity. The equations are:

(1) $$\rho_t+3H(t)\rho+\frac{1}{a(t)}\vec{\nabla}.(\rho\vec{u})=0$$

(2) $$(\rho\vec{u})_t+4H(t)\rho\vec{u}+\frac{1}{a(t)}[(\vec{\nabla}.(\rho\vec{u}))\vec{u}+\rho(\vec{u}.\vec{\nabla})\vec{u}]=-\frac{1}{a(t)}\vec{\nabla}p-\frac{\rho}{a(t)}\vec{\nabla}\Phi\ .$$

(3) $$\Delta\Phi=4\pi Ga^2\rho\ .$$

In texts of cosmology ([10] p233,[5] p 294) (2) is stated in the simpler form

$$\vec{u}_t + H(t)\vec{u} + \frac{1}{a(t)}(\vec{u}.\vec{\nabla})\vec{u} = -\frac{1}{a(t)\rho}\vec{\nabla}p - \frac{1}{a(t)}\vec{\nabla}\Phi .$$

The form (2) (obtained by multiplying (1) by $\vec{u}$, the above equation by $\rho$, and by adding the two equations so obtained) is better suited for our method and it follows directly from application of the laws of physics. We shall treat this system with the following splitting of equations, i.e. the two subsystems will be treated successively on the same cells. First subsystem:

$$\rho_t + 3H(t)\rho + \frac{1}{a(t)}\vec{\nabla}.(\rho\vec{u}) = 0$$

$$(\rho\vec{u})_t + 4H(t)\rho\vec{u} + \frac{1}{a(t)}[(\vec{\nabla}.(\rho\vec{u}))\vec{u} + \rho(\vec{u}.\vec{\nabla})\vec{u}] = 0 .$$

Second subsystem:

$$\rho_t = 0, \text{ i.e. } \rho \text{ is constant in time,}$$

$$\Delta\Phi = 4\pi G a^2 \rho ,$$

$$(\rho\vec{u})_t = -\frac{1}{a(t)}\vec{\nabla}p - \frac{\rho}{a(t)}\vec{\nabla}\Phi .$$

Note that the numbers 3 and 4 above in factor of H(t) are due to the dimension 3 of space. We shall anyway keep these numbers in the 2D tests since the 2D case has no physical significance: the 2D case is used for convenience (it works on any PC) and we only observe qualitative results.

## 2. Newtonian case, first subsystem: delta wave solution of the Riemann problem in 1D

The first system in the splitting is the following:

$$\rho_t + 3H(t)\rho + \frac{1}{a(t)}.(\rho u)_x = 0$$

(4)

$$(\rho u)_t + \ 4H(t)(\rho u) + \frac{1}{a(t)}(\rho u^2)_x = 0 .$$

The calculations of delta waves solution of the Riemann problem are difficult. On the other hand one does not need their precise formulas if one constructs the Godunov scheme from projection of delta waves by means of an extension of the method in appendix 2: formulas (6) below which are only based on the classical step function solution when it exists. So we give only numerical evidence of their existence in appendix 1.

**3. Newtonian case, first subsystem:**

Fortunately one has in a particular case a solution of the Riemann problem made of step functions so that one can apply the classical average value method for the projection of step functions in this particular case. Set:

(5) $$c_l(t)= u_l(0)a(0)\int_0^t \frac{ds}{(a(s))^2},$$

$$c_r(t)= u_r(0)\, a(0)\int_0^t \frac{ds}{(a(s))^2}.$$

If $u_l(0) < u_r(0)$ then obviously $c_l(t) < c_r(t)$ $\forall$ t>0.

**Proposition:** If $u_l(0) < u_r(0)$ one has the following solution of the Riemann problem:

if $x< c_l(t)$: $\rho(x,t)= \rho_l(0).(\frac{a(0)}{a(t)})^3$, $u(x,t) = u_l(0).(\frac{a(0)}{a(t)})$,

if $c_l(t)<x< c_r(t)$: $\rho(x,t)=0$, u(x,t) unspecified,

if $x> c_r(t)$: $\rho(x,t)= \rho_r(0).(\frac{a(0)}{a(t)})^3$, $u(x,t) = u_r(0).(\frac{a(0)}{a(t)})$.

Proof. We give the jump conditions of the discontinuity at $x= c_r(t)$. Set $\rho(x,t) = \rho_r(t)\, K(x-c_r(t))$, $u(x,t) = u_r(t)\, K(x-c_r(t))$ with K a Heaviside function (the letter H is used for the Hubble function). The first equation in (4) gives:

$\rho_r{}'K - c_r{}'\, \rho_r K' + 3H\, \rho_r K + \frac{1}{a}\, \rho_r u_r (K^2)' = 0;$

for $x> c_r(t)$ it gives $\rho_r{}' + 3H\, \rho_r=0$, hence the formula for $\rho(x,t)$. For $x= c_r(t)$ it gives $c_r{}'=\frac{1}{a}\, u_r$. The second equation in (4) gives:

$(\rho_r u_r)'\, K^2 - c_r{}'\, \rho_r\, u_r (K^2)'+4H\rho_r\, u_r\, K^2 +\frac{1}{a}\, \rho_r (u_r)^2\, (K^3)'=0;$

for $x> c_r(t)$ it gives $(\rho_r\, u_r)' + 4H\, \rho_r\, u_r=0$ hence the formula for u(x,t). For $x= c_r(t)$ it gives again $c_r{}'=\frac{1}{a}\, u_r$, hence the formula for $c_r$. One obtains similar jump conditions at $x= c_l(t)$.

Then the calculations in Appendix 2 to construct a Godunov scheme apply in a straightforward way: $c_l(\Delta t)$, $c_r(\Delta t)$ replace $u_l.\Delta t$ and $u_r.\Delta t$ respectively and use the formulas of $\rho$ and u in the above proposition to project $\rho$ and $\rho$u at time $\Delta$t=rh: here the quantities

$$\rho_l(\text{rh}) = \rho_l(0).(\frac{a(0)}{a(rh)})^3,$$

$$\rho_r(\text{rh}) = \rho_r(0).(\frac{a(0)}{a(rh)})^3,$$

$$\rho_l(\text{rh})u_l(\text{rh}) \text{ with } u_l(\text{rh}) = u_l(0).(\frac{a(0)}{a(rh)}),$$

$$\rho_r(\text{rh}), u_r(\text{rh}) \text{ with } u_r(\text{rh}) = u_r(0).(\frac{a(0)}{a(rh)}),$$

replace respectively the time independent quantities $\rho_l$, $\rho_r$, $(\rho u)_l$, $(\rho u)_r$ of Appendix 2 which correspond to $a(t) =1 \forall t$.

One obtains the analog of theorem 1 in Appendix 2 concerning physical interpretation. The formulas in remark 3 of section 4 in Appendix 2 become here (notation [a,b](x)=1 if a<x<b, 0 if not):

(6) $$\overline{\rho_l} = [\rho_l(rh)\int_0^{h/2}[-h/2, c_l(rh)](x)dx + \rho_r(rh)\int_0^{h/2}[c_r(rh), h/2](x)dx]/(h/2),$$

$$\overline{\rho_r} = [\rho_l(rh)\int_{-h/2}^{0}[-h/2, c_l(rh)](x)(x)dx + \rho_r(rh)\int_{-h/2}^{0}[c_r(rh), h/2](x)dx]/(h/2),$$

same for $\overline{(\rho u)}_l$ and $\overline{(\rho u)}_r$ using $\rho_l(\text{rh}).u_l(\text{rh})$ and $\rho_r(\text{rh}).u_r(\text{rh})$ in place of $\rho_l(\text{rh})$ and $\rho_r(\text{rh})$ respectively. Let us recall from Appendix 2 that these formulas have been obtained in case $c_l(t)< c_r(t)$ from the classical Godunov method and that our method consists in adopting them in the case $c_l(t)>c_r(t)$ as their "natural" extension.

## 4. Newtonian system: the numerical scheme in any space dimension

For the first subsystem it is an adaptation of the one in Appendix 2, taking into account the above formulas due to expansion, see appendices 1 and 2. We set $t_n = n\Delta t$. From (6), (8) and the above proposition:

**Formulas of the 1D scheme.** The mean values $\rho_i^{n+1}$ and $(\rho u)_i^{n+1}$ in the cell [i-h/2, i+h/2] at time $t_{n+1}$ are given by the following formulas that sum up the possible contributions of the cells i-1,i and i+1. At first auxiliary formulas:

$$\overline{\rho_i^n} = \rho_i^n (\frac{a(t_n)}{a(t_{n+1})})^3, \quad \overline{(\rho u)_i^n} = (\rho u)_i^n (\frac{a(t_n)}{a(t_{n+1})})^4,$$

$$\overline{u_i^n}\,\text{rh} = u_i^n a(t_n) \int_{nrh}^{(n+1)rh} \frac{ds}{(a(s))^2}.$$

For this integral we use the trapeze formula for integration:

$$\overline{u_i^n} = u_i^n a(t_n)\frac{1}{2}\left(\frac{1}{a(t_{n+1})^2}+\frac{1}{a(t_n)^2}\right).$$

Let us recall from Appendix 2 that

L(a,b) = length of [0,1] $\cap$ [a,b] = min (1,b) –max (0,a).

Then the scheme is (see Appendix 2):

$$\rho_i^{n+1} = \overline{\rho}_{i-1}^{\,n}\, L(-1+r\overline{u_{i-1}^n}, r\overline{u_{i-1}^n},) + \overline{\rho}_i^{\,n}\, L(r\overline{u_i^n}, 1+r\overline{u_i^n}) + \overline{\rho}_{i+1}^{\,n}\, L(1+r\overline{u_{i+1}^n}, 2+r\overline{u_{i+1}^n}), \quad (10)$$

$$(\rho u)_i^{n+1} = \overline{\rho u_{i-1}^n}\, L(-1+r\overline{u_{i-1}^n}, r\overline{u_{i-1}^n},) + \overline{\rho u_i^n}\, L(r\overline{u_i^n}, 1+r\overline{u_i^n}) + \overline{\rho u_{i+1}^n}\, L(1+r\overline{u_{i+1}^n}, 2+r\overline{u_{i+1}^n}),$$

completed by

$$u_i^{n+1} = \frac{(\rho u)_i^{n+1}}{\rho_i^{n+1}}.$$

The CFL condition is

$$r \max\{|u_i^n|\} \le 1.$$

For the second system in the splitting the energy density $\rho$ is constant in time, so the formulas for $\Phi$, then $\rho$u, are explicit in case of a state law p=constant.$\rho$.

The 2D and 3D cases are similar as explained in Appendix 2.

**Validation:** First one checked that in absence of expansion the scheme reduces to the one in Appendix 2. Then one checked that in absence of transfer through the boundary of the window the quantities $a(t)^3 \iint \rho(x,t)dx$, $a(t)^4 \iint \rho(x,t)u(x,t)dx$ and $a(t)^4 \iint \rho(x,t)v(x,t)dx$ are independent on time (which follows at once from the equations). This is done by imposing a density=0 near the boundary of the window. The concrete physical interpretation of the numerical scheme for (4) (same as theorem1 in Appendix 2: free streaming at the interfaces of meshes then mean values on each cell) is also a positive indication. One has also compared the result from this 1D scheme with results from the viscosity scheme (this last one is solution with second members $\varepsilon\rho_{xx}$ and $\varepsilon.(\rho u)_{xx}$ respectively with $\varepsilon > 0$ as small as possible), in cases the viscosity scheme gives a clear enough result.

## 5. Newtonian system: numerical tests in a 2D toy model

If the length scale of initial perturbations of dark matter or baryonic matter after decoupling is smaller than the effective cosmological horizon a Newtonian treatment is expected to be valid. First we consider the case of a pressureless fluid such as dark matter. Dark matter is made of weakly interacting particles which do not couple even to radiation: they do not feel photon pressure. For slow enough expansion, or weak enough pressure, one observes formation of the familiar network of cluster-filaments-voids, figures 1 and 2. All figures below are given in comoving coordinates, i.e. they represent regions with very different physical sizes (expanding with time) that had the same size at time=0.

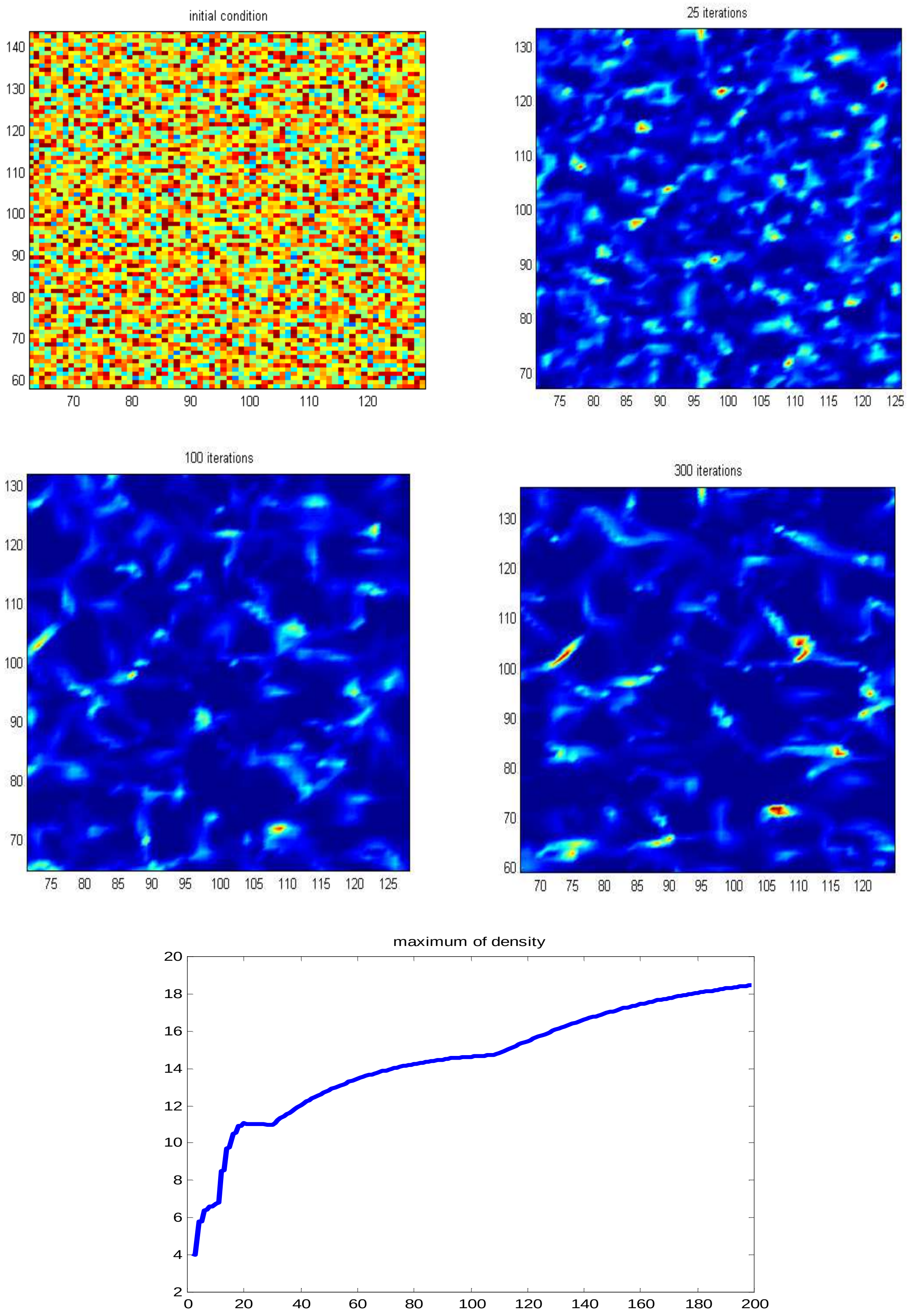

**figure 1. A film of structure formation in a slowly expanding universe made only of Newtonian matter (no pressure). Random initial conditions. The scale factors are respectively 1, 1.12, 1.5 and 2.5. The development of a filament-cluster-void network is clearly seen. Bottom: growth of the maximum value of density with time: peaks of matter grow very fast at the beginning, then, when they are separated by void regions, their growth is slow.**

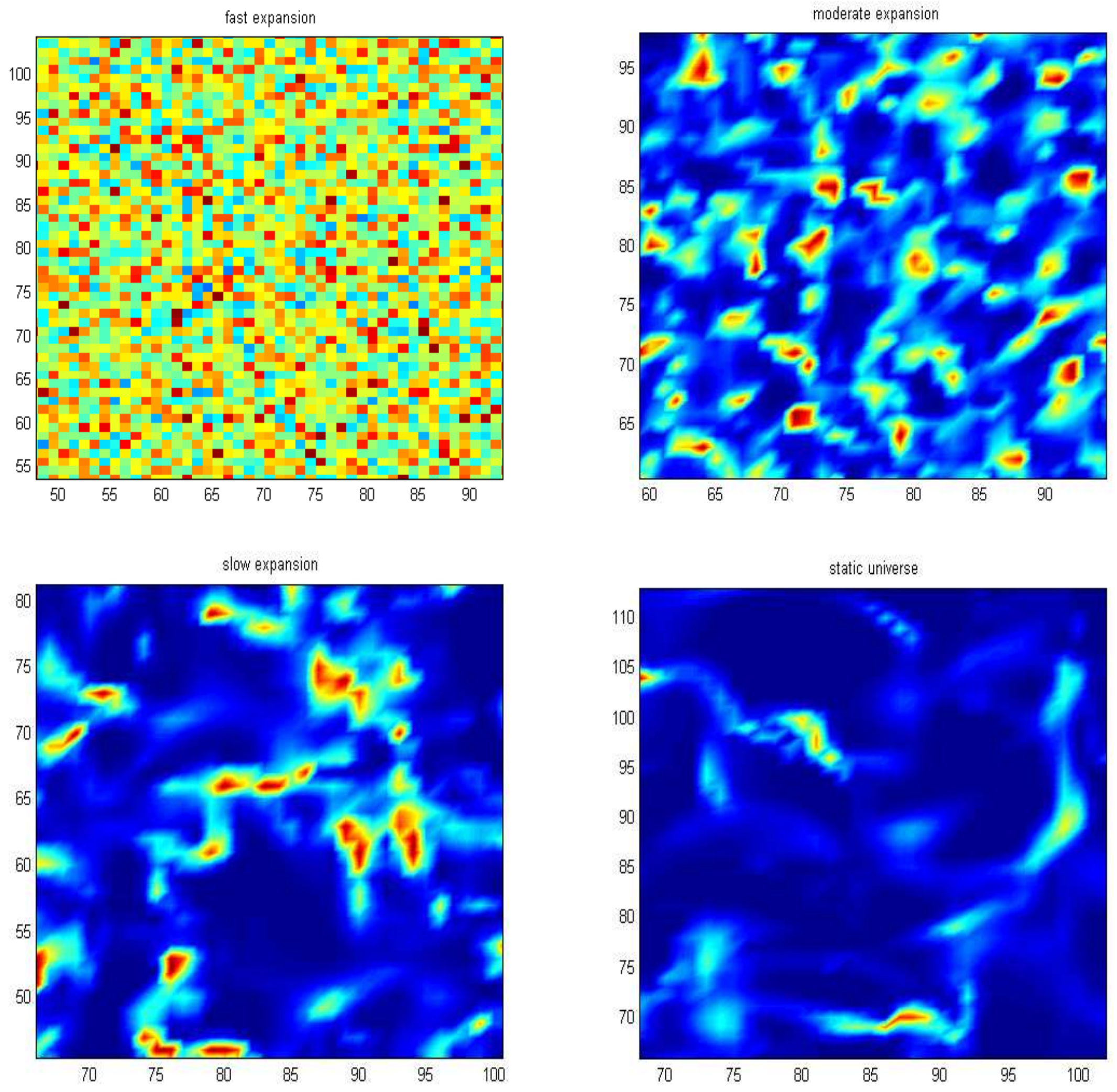

**figure 2. Newtonian system. How structure formation depends on expansion (no pressure): 100 iterations ( from top-left to bottom-right the scale factor is 128, 13.7, 3.5, 1 respectively). Random initial conditions. Structure formation is made impossible by too fast expansion, see also figure 3 where an initial structure is frozen by expansion. In case of slow enough expansion one observes the familiar network of clusters and filaments surrounding large void regions as in absence of expansion.**

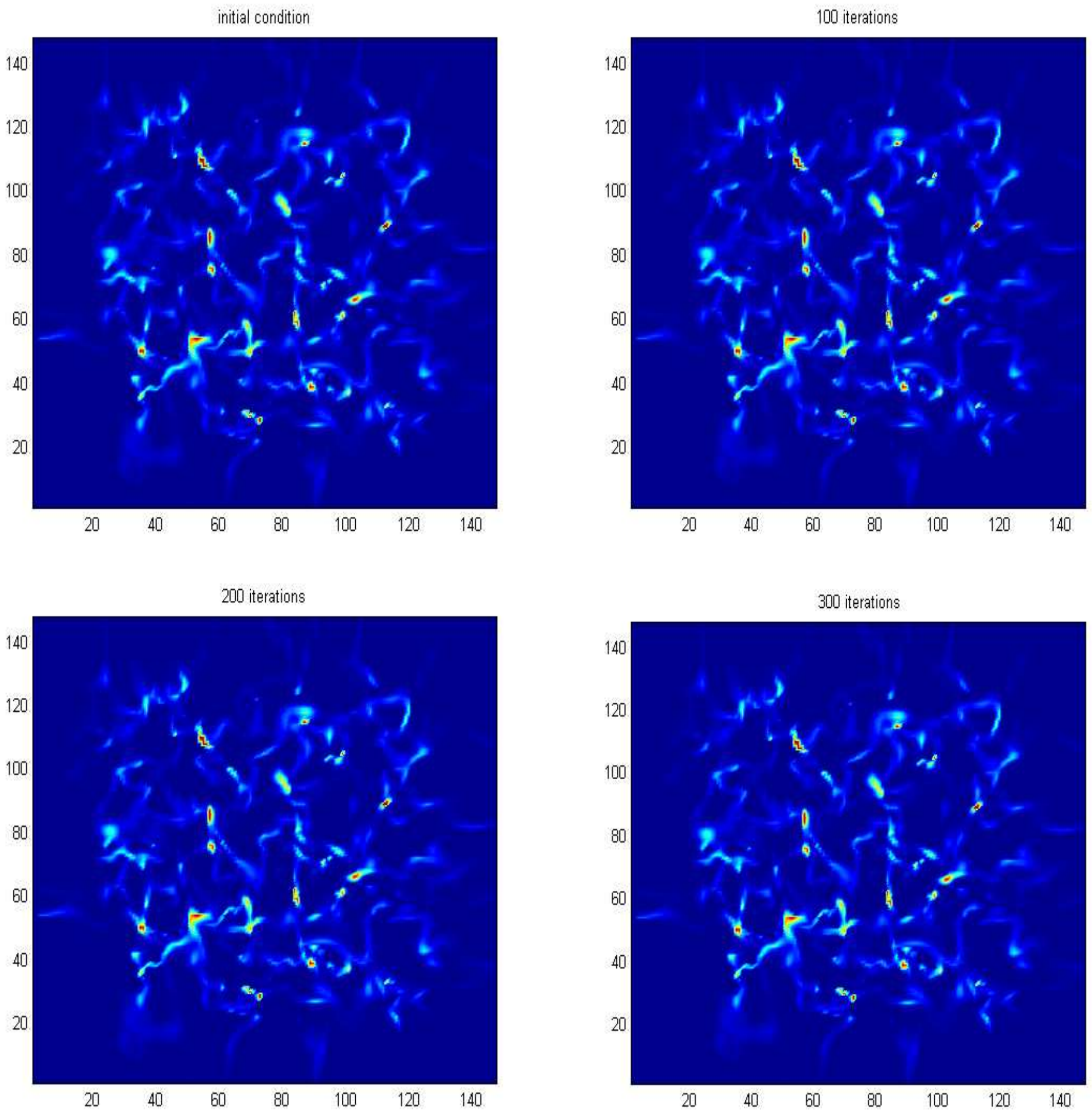

**figure 3. A structure formation is frozen by expansion. We choose as initial condition in the Newtonian system a structure that has been previously created (in a slow enough expansion background). If the rate of expansion is fast enough one observes that the structure is frozen. This is the 2D analog of the “Meszaros effect” or “stagnation”: fluctuations of Newtonian matter that enter the horizon before equivalence are frozen till equivalence. Radiation domination ensures that the universe expands very quickly (from the Friedman equations) and so Newtonian perturbations that enter through the horizon are frozen as shown in this test. In this test radiation is considered as a smooth unclustered background that affects only the overall expansion rate. In figure 9 below this test will be reproduced in presence of radiation as this is really the case. Other tests using the relativistic system show that relativistic fluctuations are also frozen: the “Meszaros effect” has been observed whether one uses the Newtonian or the relativistic description of dark matter before equivalence.**

Now we consider a Newtonian fluid having a pressure of the form p=const.$\rho$ as a simplification of the formula $p = \frac{1}{3} < v^2 > \rho$.

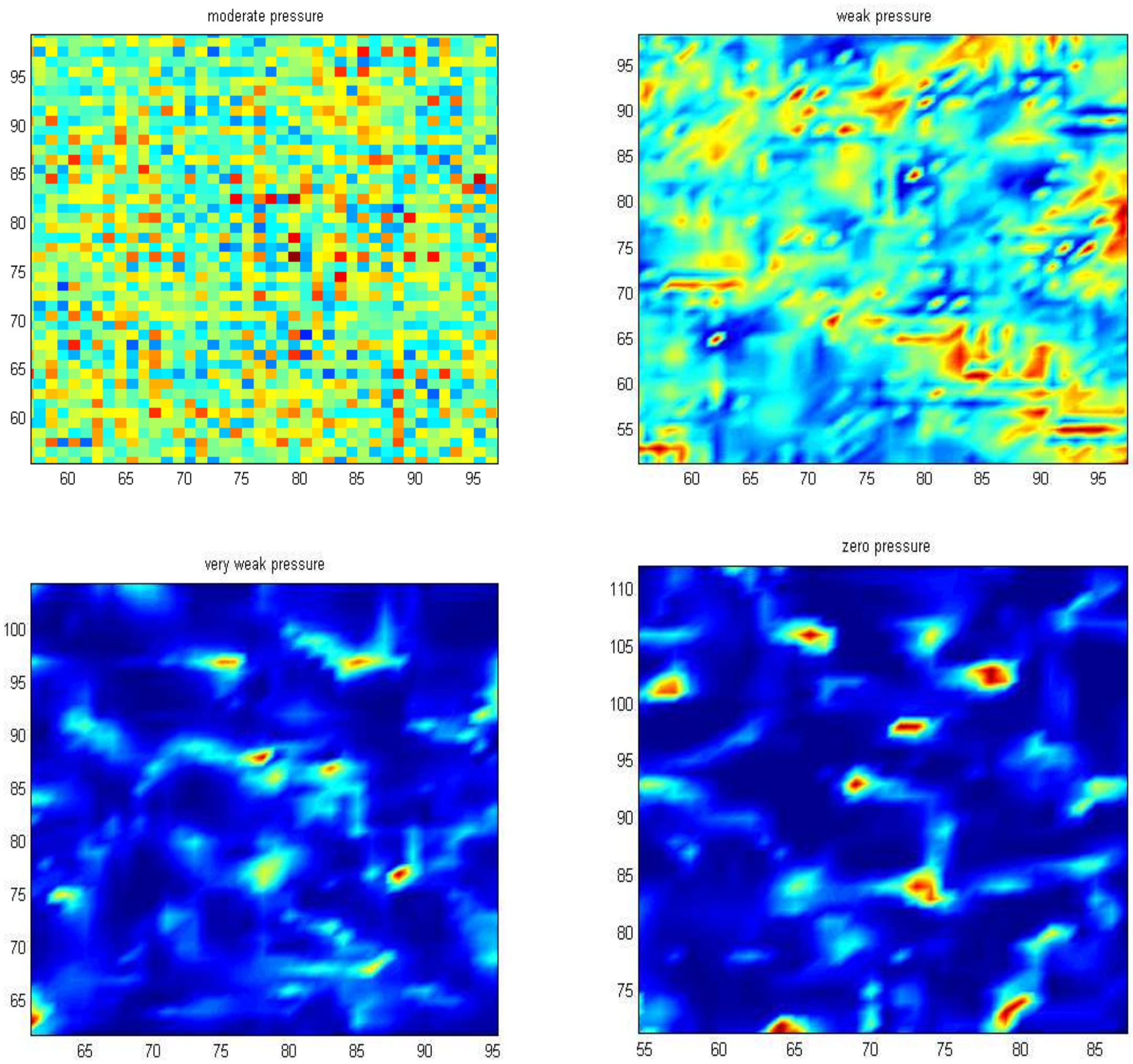

**figure 4. Newtonian system. Influence of pressure on structure formation. These pictures are given in a case of slow expansion (multiplication of the scale factor by 3.5 in 100 iterations). One observes that structure formation requires a weak enough pressure (state law p=const.$\rho$ where from left to right and top to bottom const= 10, 1, 0.1, 0 respectively). One knows that the Jeans' length provides an information on the gravitational collapse: perturbations of constant density collapse iff their size is > Jeans' length; in the initial condition the values of $\rho$ are at random between 0.9 and 1.1, therefore approximately 1 and from left to right and top to bottom the Jeans' length is approximately 0.5, 0.05, 0.005, 0.The numerical results are in agreement with the fact structure formation depends crucially on the Jeans'length. Indeed usual perturbations in density involving 4 meshes forming a square correspond to a diameter of 0.01 approximately; they collapse in the bottom figures. Larger rare perturbations ( larger than a square made of 49 meshes) are requested in top-right (in order to collapse). Collapse is practically impossible in top-left: indeed it is not observed.**

## 6. The equations in the relativistic domain

From the Einstein equations one obtains the set of equations ([5] p 221+ [10] p233):

(11) $$\rho_t + 3H\rho + \frac{1}{a}\nabla.[(\rho + \frac{p}{c^2})\vec{u}] = 0 ,$$

(12) $$(\rho + \frac{p}{c^2})(\vec{u}_t + H\vec{u} + \frac{1}{a}(\vec{u}.\nabla)\vec{u}) + \frac{1}{a}\nabla p + \frac{1}{a}(\rho + \frac{p}{c^2})\nabla\Phi = 0 ,$$

(13) $$\Delta\Phi = 4\pi G a^2 (\rho + 3\frac{p}{c^2}) .$$

In the sequel we mainly consider the state law $p = \frac{1}{3}c^2\rho$ of radiation. Then easy calculations give the system of scalar equations:,

(14) $$\rho_t + 3H\rho + \frac{4}{3a}[(\rho u)_x + (\rho v)_y + (\rho w)_z] = 0 ,$$

(14x)

$$(\rho u)_t + 4H\rho u + \frac{4}{3a}[(\rho u^2)_x + (\rho uv)_y + (\rho uw)_z] - \frac{1}{3a}\rho(uu_x + vu_y + wu_z) + \frac{c^2}{4a}\rho_x + \frac{1}{a}\rho\Phi_x = 0$$

(14y)

$$(\rho v)_t + 4H\rho v + \frac{4}{3a}[(\rho uv)_x + (\rho v^2)_y + (\rho vw)_z] - \frac{1}{3a}\rho(uv_x + vv_y + wv_z) + \frac{c^2}{4a}\rho_y + \frac{1}{a}\rho\Phi_y = 0$$

(14z)

$$(\rho w)_t + 4H\rho w + \frac{4}{3a}[(\rho uw)_x + (\rho vw)_y + (\rho w^2)_z] - \frac{1}{3a}\rho(uw_x + vw_y + ww_z) + \frac{c^2}{4a}\rho_z + \frac{1}{a}\rho\Phi_z = 0$$

(15) $$\Delta\Phi = 8\pi G a^2 \rho .$$

## 7. Relativistic system in 2D: splitting

We are forced to use a splitting into 3 subsystems. The first term in the splitting is the one in section 3 in which the velocity has been multiplied by $\frac{4}{3}$:

$$\rho_t + 3H\rho + \frac{4}{3a}[(\rho u)_x + (\rho v)_y] = 0$$

$$(\rho u)_t + 4H\rho u + \frac{4}{3a}[(\rho u^2)_x + (\rho uv)_y] = 0$$

$$(\rho v)_t + 4H\rho v + \frac{4}{3a}[(\rho uv)_x + (\rho v^2)_y] = 0 .$$

The second system in the splitting is:

$\rho$ =constant in time and $\Delta\Phi = 8\pi G a^2 \rho$

$$u_t = -\frac{c^2}{4a}\frac{\rho_x}{\rho} - \frac{1}{a}\Phi_x$$
$$v_t = -\frac{c^2}{4a}\frac{\rho_y}{\rho} - \frac{1}{a}\Phi_y \, .$$

It is solved similarly as in the Newtonian case: standard solution of the Poisson equation and integration in time for u and v. The third system in the splitting is:

$$u_t = \frac{1}{3a}(uu_x + vu_y)$$
$$v_t = \frac{1}{3a}(uv_x + vv_y) \, .$$

It is solved according to the method used in Appendix 2 for pressureless fluid dynamics, considering here a transportation of each cell by the vector $(-\frac{u}{2},-\frac{v}{2})\frac{\Delta t}{3a}$, which follows from these equations.

The numerical tests show that this scheme usually demands a smaller value of r (to fulfill the two CFL conditions in step 1 and step 3) than the scheme for the Newtonian system, and that this value diminishes strongly when the value of c increases.

**8. Relativistic system: numerical tests in a 2D toy model**

In modern cosmological theories involving inflation the relativistic treatment is extremely important since fluctuations outside the horizon must be handled using general relativity [5, §10.12]. They are supposed to play an important role in structure formation by crossing the horizon, being frozen till equivalence and starting growth after equivalence (we are going to observe all these steps from our 2D toy models).

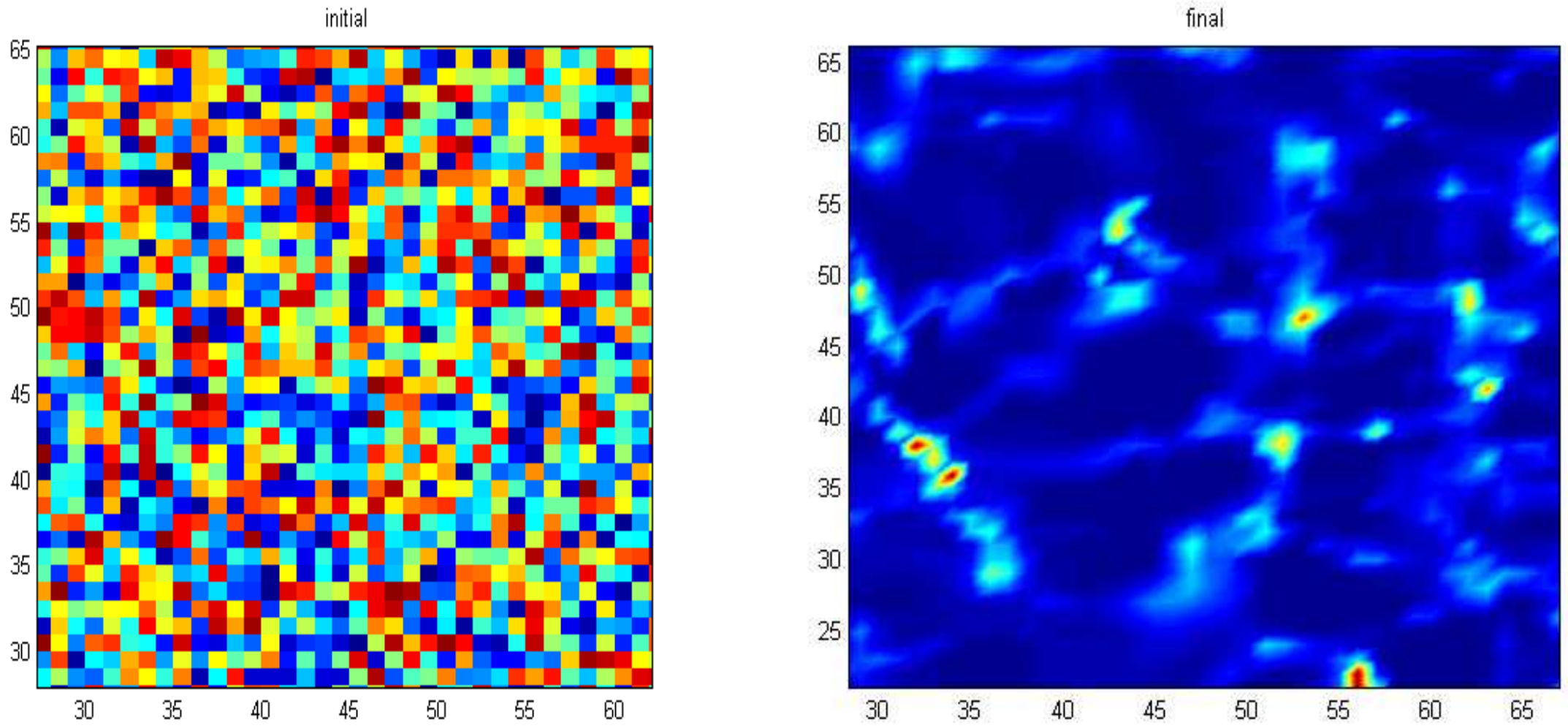

**figure 5 . Formation of structures from the 2D-relativistic system in case of slow expansion and small value of the velocity of light (that can be interpreted as a model for super-horizon perturbations; these perturbations have to be handled using the relativistic system and the small values of c correspond to very large initial perturbations). The initial fluid is at random; state law of radiation** $p = \frac{1}{3} c^2 \rho$ **. This suggests that super-horizon structures (in radiation, baryonic matter and dark matter) could be created, then enter into the horizon when the universe expands. One has observed these structures are frozen by fast expansion exactly like the Newtonian ones in figure 3 (Meszaros effect).**

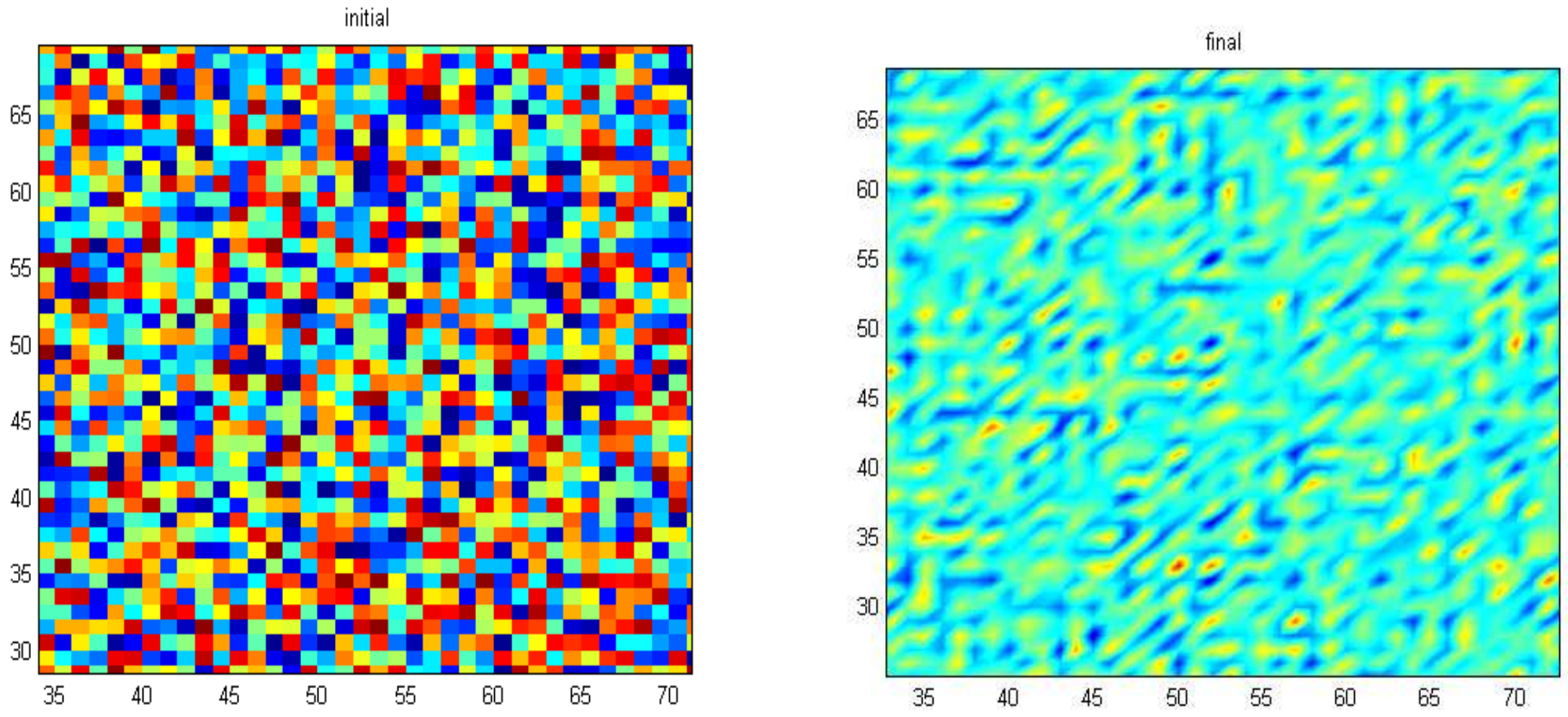

**figure 6. The relativistic system does not form structures when the velocity of light is not small (case of sub-horizon perturbations of radiation at all epochs , and of baryonic matter before decoupling) even in absence of expansion (the visualization on the right is the one of matter at random, as on the left; it differs only by the presentation which is the one used to show structures).**

## 9. Multifluid tests

The numerical schemes for each fluid are only superposed independently except the gravitation step that works with a density of energy which is the sum of the densities of the two fluids in the Poisson equation. We do not give them in appendix since they can be at once constructed from the single fluid schemes.

One can consider the coupled (through gravitation) evolution of densities of dark matter (Newtonian) and baryons (relativistic) before decoupling. Baryons are then coupled to photons. Dark matter perturbations are growing in this epoch and exert a gravitational pull of baryons. Since baryons are tightly coupled to photons they have a huge internal pressure which opposes to dark matter pull.

Below we consider the evolution of perturbations of a non-relativistic component (dark matter) in a universe containing also a relativistic component (radiation and baryons coupled to radiation): the non-relativistic component forms structures from a random initial distribution while the relativistic component remains unperturbed (its distribution remains at random), [5,§10.11].

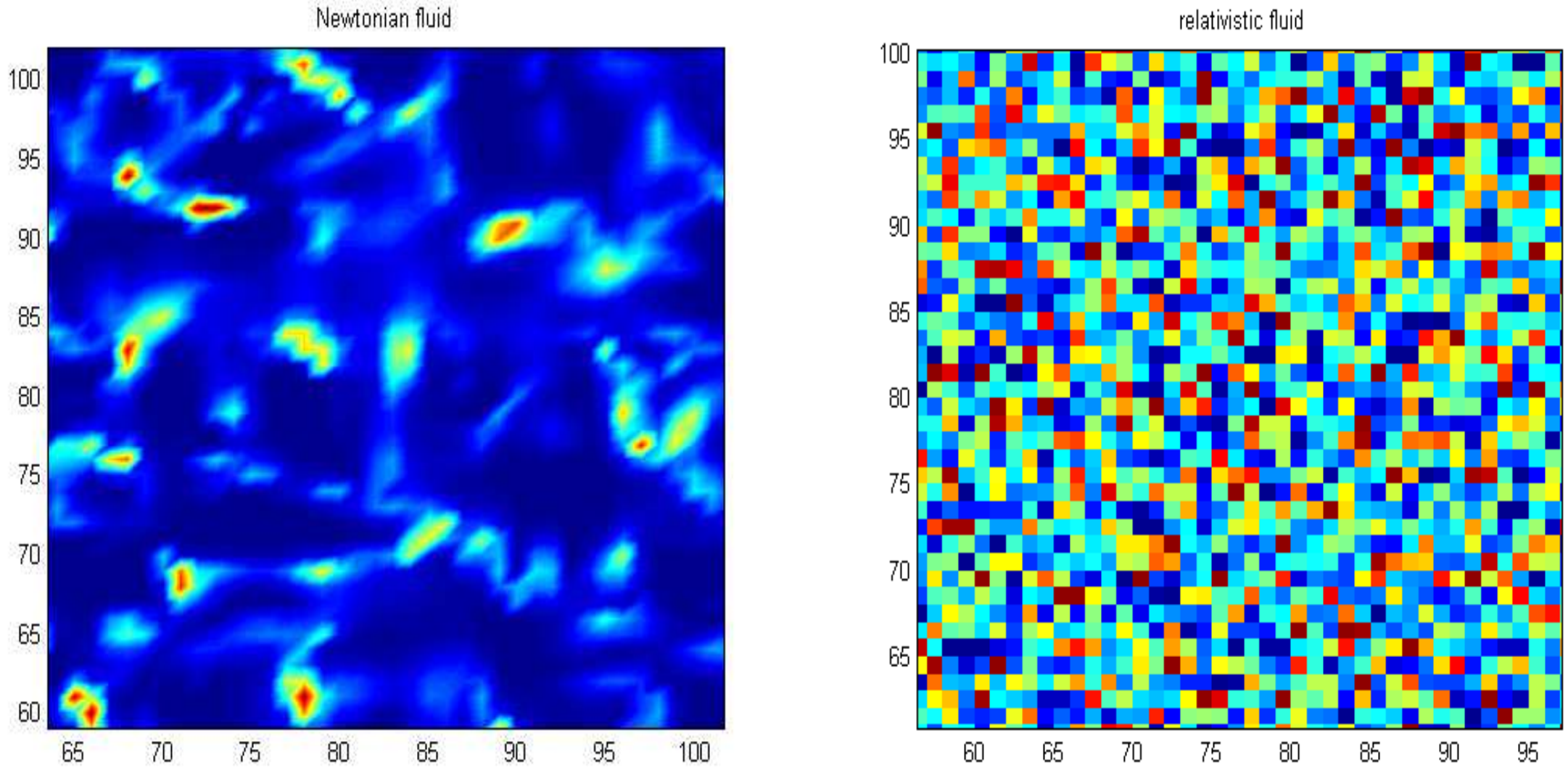

**figure 7 . Mixture of a pressureless Newtonian fluid and of a relativistic fluid (50% each, epoch of the equivalence); random initial data for both fluids; slow expansion: scale factor 3.5 in 100 iterations. The presence of relativistic fluid has not prevented structure formation in the Newtonian fluid; this suggests structure formation in dark matter has developed between equivalence and decoupling, so that a structure of dark matter (not directly observable) would exist at the time of the CMB, as well as the (observed) near absence of structure of baryonic matter (relativistic before decoupling). See the behaviour of this mixture after the time of the CMB in figures 9, 10 below.**

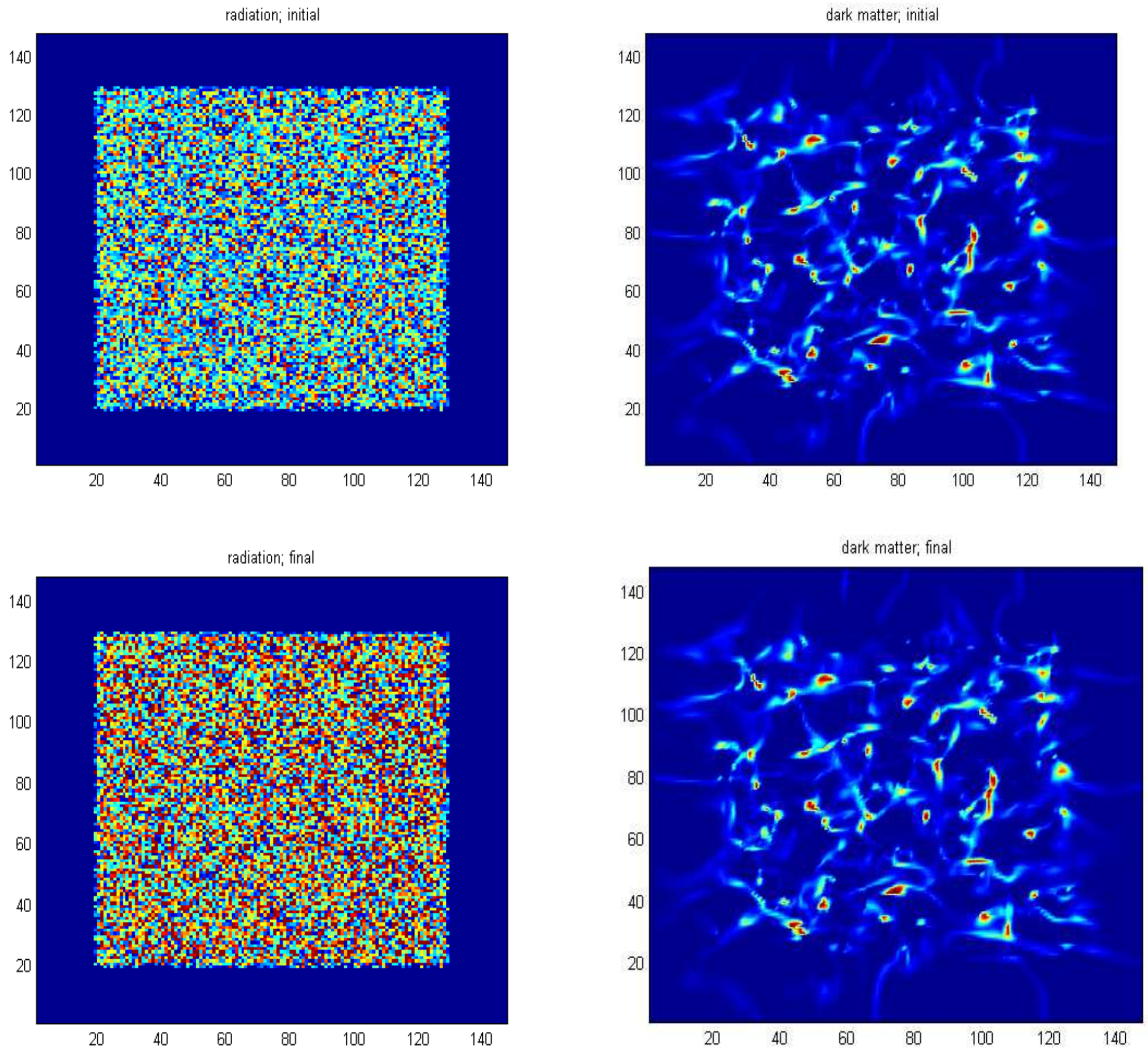

**figure 8. A simulation of the Meszaros effect more complete than the one in figure 3. Here radiation is taken into account using the mixed Newtonian-relativistic scheme: 90% of radiation and 10 % of dark matter, fast expansion. One observes that the dark matter structure is frozen due to the fast inflation. One deduces from this test that the presence of radiation plays no direct role (but an indirect role: the dominant energy of radiation drives the universe to expand fast from the Friedman equations). Radiation remains at random.**

At the epoch of decoupling, baryons have not yet formed structures while structures of dark matter already exist and serve as gravitational potential wells for the emergence of structures of baryonic matter.

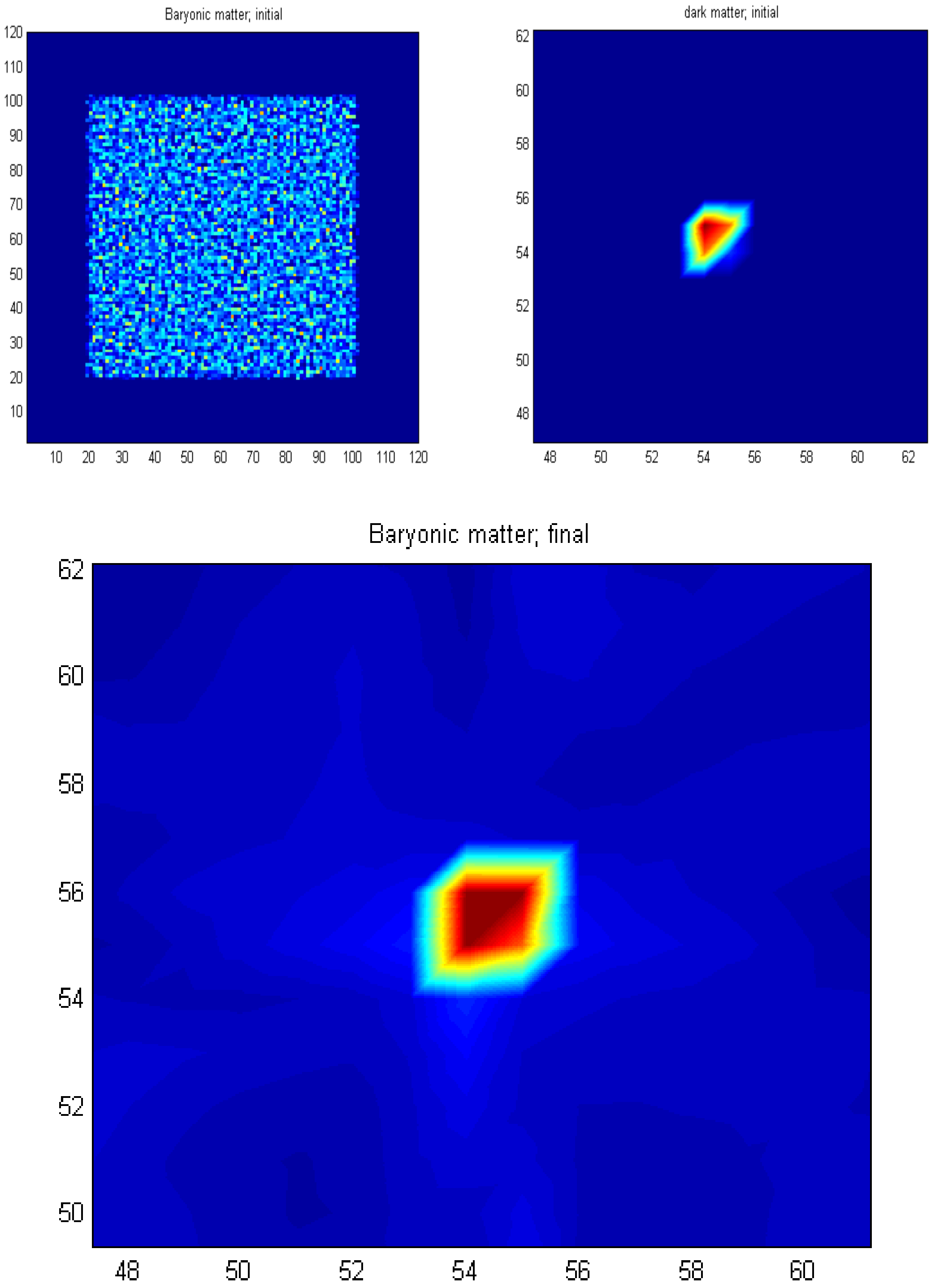

**figure 9. 2D Mixture of 2 Newtonian fluids: dark matter (80%) and baryonic matter (20%). Initial condition: top-left: baryonic matter at random in a large region and top-right: a peak of dark matter previously formed. Final result after a few time steps: one observes that most of baryonic matter is agglomerated on the peak of dark matter outside of which there subsists a diluted "gas" of baryonic matter. For a clearer visualization in 1D see figure below.**

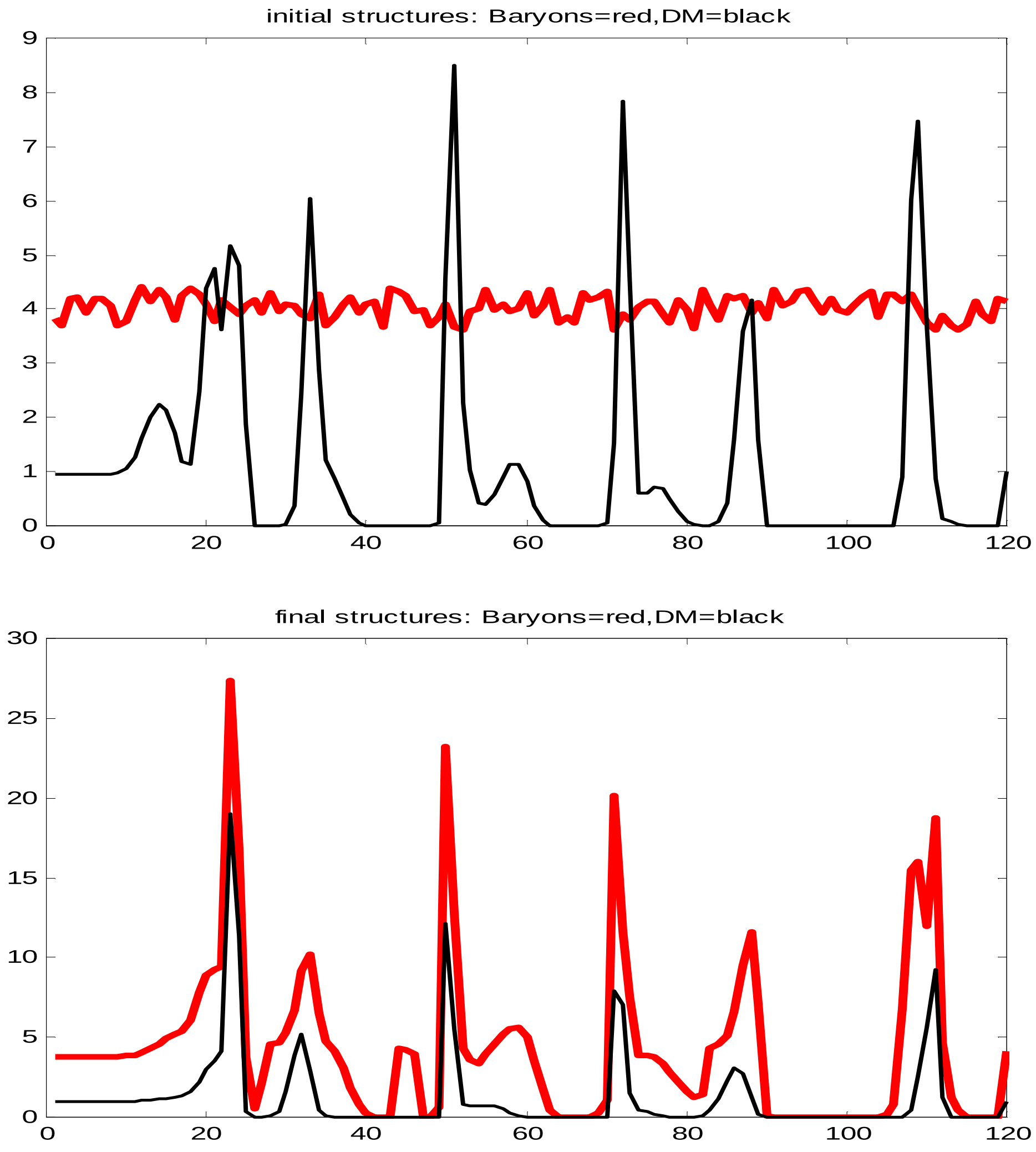

**figure 10. 1D Mixture of 2 Newtonian fluids: dark matter (80%) and baryonic matter (20%) (the vertical scale of baryonic matter has been enlarged for a better visualization). Initial condition: baryonic matter at random and structures of dark matter previously formed (with their own velocity). After a few iterations one observes the creation of baryonic structures that mimick perfectly well the structures of dark matter.**

## 10. Taking entropy into account

Entropy plays a significant role in some cosmological models [5 p 207,216]. In the perfect fluid approximation (which means neglecting viscous forces in the fluid) the Euler equations in a static background (for simplification) are stated as:

$$\rho_t + \vec{\nabla}.(\rho \vec{u}) = 0$$

$$\vec{u}_t + (\vec{u}.\vec{\nabla})\vec{u} = -\frac{1}{\rho}\vec{\nabla}p - \vec{\nabla}\Phi$$

$$\Delta\Phi = 4\pi G a^2 \rho$$
$$p = p(\rho, S)$$
$$S_t + \vec{u}.\vec{\nabla} S = 0.$$

where S is the entropy per unit mass: the last equation is “theorem H” of Boltzmann for a collisionless system or a system in thermodynamical equilibrium.

In perturbation theory one uses a state law $dp = \gamma^2 d\rho + \alpha^2 \frac{dS}{S}$ which suggests the choice of a state law $p = \gamma^2 \rho + \alpha^2 \log(S)$. For a Godunov type scheme one has to consider shock waves: then the term $\vec{u}.\vec{\nabla} S$ does not make sense in the theory of distributions. A mathematical study shows that the entropy equation could be treated by the method of Appendix 2 in sections 3, 4.

## 11. Intense stellar formation inside a shock wave in a dust cloud

Shock waves can collect matter (the delta peaks in 1D tests). On the surface of the shock wave the situation looks like the situation in a 2D contracting universe. In this case, one observes from the Newtonian system intense structure formation, even in presence of significative pressure that would inhibit structure formation in a static universe.

## Conclusion

Recapitulation of the numerical observations from these 1D and 2D toy models (there is no insurance one has performed all tests in the most adequate conditions: it is not convenient on a PC to give their physical values to the gravity constant and to the velocity of light, and impossible to perform efficient 3D tests).

**Super-horizon perturbations**, figure 5:
formation of structures from the relativistic system (radiation, baryonic matter and dark matter). They can enter into the horizon from the fast expansion before equivalence.

**Sub-horizon perturbations before equivalence,** figures 2, 3, 6, 8:
structures of dark matter that cross the horizon are frozen (Meszaros effect, observed with the relativistic and also the Newtonian system); no structure formation in radiation, baryonic matter (relativistic system) and dark matter ( Newtonian system but frozen).

**Sub-horizon perturbations between equivalence and decoupling,** figures 2, 4, 6, 7: structures of dark matter grow (Newtonian system); no structure formation in baryonic matter and radiation (relativistic system). The role of pressure (Jeans' theory) has been observed in the Newtonian system (and also through non small values of c in the relativistic system).

**Sub-horizon perturbations after decoupling,** figures 2, 6, 9, 10:
structures of dark matter grow; baryonic matter catch up with the previously existing structures of dark matter. No structure formation in radiation.

Therefore from the full nonlinear systems of fluid dynamics for cosmic fluids one recovers classical results of the theory of cosmological perturbations.

It seems possible that (in hands of specialists in cosmology) these schemes could be at the origin of a simple and efficient "machine" for handling nonlinear evolution of cosmological fluids capable to predict the distribution of galaxies and other structures by testing various hypotheses (in the Cold Dark Matter model for instance, see [5 §15.7]. Their 3D extension is immediate without dimensional splitting (i.e. without loss in efficiency).

## Appendix 1: Delta waves solutions of the Riemann problem for pressureless fluid dynamics in an expanding background

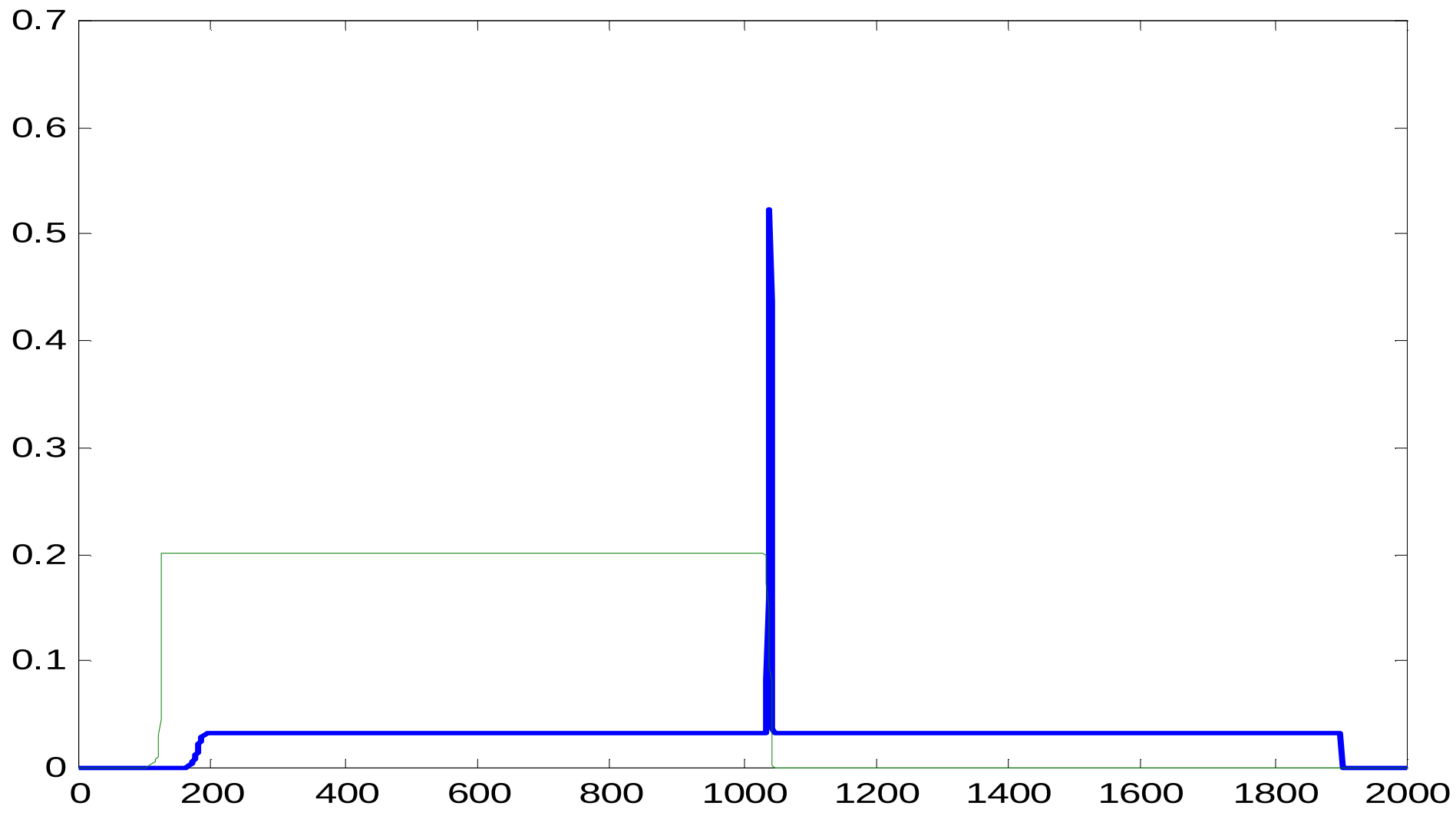

**<u>figure 12:</u> One observes delta peaks in the numerical solution of system (4) in a case of compression: blue=density, green=velocity**

## Appendix 2. A Method of projection of delta waves in a Godunov scheme for pressureless fluid dynamics

In [1] the authors noticed that the solution of the Riemann problem for the system of pressureless fluid dynamics

$$\text{(S)} \qquad \begin{aligned} &\rho_t + (\rho u)_x = 0, \\ &(\rho u)_t + (\rho u^2)_x = 0, \end{aligned}$$

shows a delta wave located on the discontinuity of the solution. Nevertheless, they succeeded to extend the Godunov method to this case, and obtained excellent numerical results. After the pioneering article [1], various numerical methods have been proposed for the numerical solution of system (S). References are given in [2,3,4,7]. As far as the A. knows they are different from the method presented here, which appears as a modification of the one in [1] in case of numerous changes in sign of velocity (in the applications to large structure formation in cosmology the sign of initial velocity in each cell varies at random).

In the projection step of the Godunov method one has to share the delta waves occurring in the Riemann problems into left-hand-side and right-hand-side contributions. In the solution of the Riemann problem there occur two cases. In the first case, we have only one solution of the Riemann problem. It is made of a physically meaningful delta wave, that we do not know *a priori* how to share. In the second case, one has two possible solutions: a physical one that has a classical form (step functions without delta wave), and a non-physical one involving a delta wave. In this second case, one obtains a Godunov scheme from the physical solution in form of step functions, which permits to compute the formulas governing the sharing of the non-physical delta wave that would lead to the same scheme. The method in this paper consists in applying the same formulas in the first case, when the unique solution is in form of a delta wave.

The scheme so obtained coincides with the scheme of [1] in the cases they considered in their applications. Stability of the scheme is proved as well as a global convergence result in any dimension and for any configuration of the velocity field. The scheme in this paper has been adapted in [6] for numerical approximations of standard systems of cosmology modelling Newtonian and relativistic fluids coupled to gravitation in an expanding background, providing numerical simulations of large structure formation in cosmology.

### 1. Description of the numerical scheme in [1]

Standard 1D notation is used: the space cells are the segments [*ih-h/2, ih+h/2*], $i\in Z$, the space step is denoted by *h* and the time step by $\Delta t$; we set $t_n = n.\Delta t$ and $r = \Delta t/h$. The constant values of $\rho$ and $u$ on the cell [*ih-h/2, ih+h/2*] at time $t_n$ are denoted by $\rho_i^n$ and $u_i^n$. In the scheme in [1] the passage from $\{\rho_i^n,\ u_i^n\}$ to $\{\rho_i^{n+1},\ u_i^{n+1}\}$ is done as follows. One introduces three intermediate values (attached to the junctions of cells)

$$w_{i+1/2}^n = \sqrt{\rho_i^n}.u_i^n + \sqrt{\rho_{i+1}^n}.u_{i+1}^n$$

(1) and

$$u_{i+1/2}^n,\ \ \rho_{i+1/2}^n,$$

defined by:

- if $u_i^n>0$ and $u_{i+1}^n>0$ then $u_{i+1/2}^n = u_i^n$, $\rho_{i+1/2}^n = \rho_i^n$,
- if $u_i^n>0$ and $u_{i+1}^n<0$
  - if $w_{i+1/2}^n>0$ then $u_{i+1/2}^n = u_i^n$, $\rho_{i+1/2}^n = \rho_i^n$,
  - if $w_{i+1/2}^n<0$ then $u_{i+1/2}^n = u_{i+1}^n$, $\rho_{i+1/2}^n = \rho_{i+1}^n$,
- if $u_i^n<0$ and $u_{i+1}^n>0$ then $u_{i+1/2}^n = 0$, $\rho_{i+1/2}^n = 0$,
- if $u_i^n<0$ and $u_{i+1}^n<0$ then $u_{i+1/2}^n = u_{i+1}^n$, $\rho_{i+1/2}^n = \rho_{i+1}^n$.

Finally one computes the values $\{\rho_i^{n+1}, u_i^{n+1}\}$ from the formulas

$$\begin{aligned} &\rho_i^{n+1} = \rho_i^n - r\rho_{i+1/2}^n\, u_{i+1/2}^n + r\rho_{i-1/2}^n u_{i-1/2}^n, \\ (2)\qquad &(\rho u)_i^{n+1} = \rho_i^n\, u_i^n - r\rho_{i+1/2}^n\,(u_{i+1/2}^n)^2 + r\rho_{i-1/2}^n (u_{i-1/2}^n)^2, \\ &u_i^{n+1} = (\rho u)_i^{n+1}/\rho_i^{n+1}. \end{aligned}$$

In the case $u_i^n>0$ and $u_{i+1}^n<0$, if $w_{i+1/2}^n=0$, then the two possible values of $\rho_i^{n+1}$ differ by a quantity $r(\rho_i^n|u_i^n| + \rho_{i+1}^n|u_{i+1}^n|)$. As an example let the values of $\rho_{i-1}^n$, $\rho_i^n$, $\rho_{i+1}^n$ be equal ($=\rho$), $u_{i-1}^n= 1$, $u_i^n=1$ and $u_{i+1}^n= -1-\varepsilon$; then one computes that the scheme gives $\rho_i^{n+1}= (1+2r+r\varepsilon)\rho$. Now, if one changes $u_i^n$ and $u_{i+1}^n$ into $u_i^n= 1+\varepsilon$ and $u_{i+1}^n= -1$, one computes $\rho_i^{n+1}= (1-r\varepsilon)\rho$, which differs from the previous value by a quantity $2r\rho$ (when $\varepsilon\to 0$, which makes the two possibilities undistinguishable while $2r\rho$ is not at all small). Similar trouble can be found in other schemes for system (S).

This does not cause any trouble in [1], since in the applications under consideration, changes in sign of velocity are rare events, so that the above cannot influence the final result. But, if one wants to model structure formation in cosmology, the sign of velocity changes at random in the initial conditions, thus displaying the above trouble as the general case. This motivates the search of a numerical scheme that would work also when velocity changes sign freely.

## 2. Solution of the Riemann problem

The calculations of the solution of the Riemann problem for the system of pressureless fluid dynamics [1, 3, 7] are recalled in appendix 3. The values of ($\rho$, $u$) are ($\rho_l$, $u_l$) on the left-hand-side of the discontinuity located at $x$=0 and ($\rho_r$, $u_r$) on the right-hand-side. If $w$ is any variable, we set $\Delta w = w_r - w_l$. We set

$$\begin{aligned} &u(x,t) = u_l + \Delta u.\, H(x-ct), \\ &\rho(x,t) = \rho_l + \Delta\rho .H(x-ct) + \alpha t \delta(x-ct), \\ &(\rho u)(x,t) = (\rho u)_l + \Delta(\rho u).H(x-ct) + \beta t \delta(x-ct), \\ &(\rho u)_l = \rho_l u_l,\ (\rho u)_r = \rho_r u_r,\ u=(\rho u)/\rho, \end{aligned} \tag{3}$$

where $H$ is the Heaviside step function, $\delta$ is the Dirac delta function. The velocity $u$ is discontinuous at $x=ct$, while $\rho$ and $\rho u$ display a delta peak on the discontinuity, which is proportional to time.

Calculations recalled in appendix 3 give:

$$\begin{aligned} &c = (\sqrt{\rho_r}.u_r + \sqrt{\rho_l}.u_l)/(\sqrt{\rho_r} + \sqrt{\rho_l}), \\ &\alpha = -\sqrt{\rho_l \rho_r}.\Delta u, \\ &\beta = c\alpha . \end{aligned} \tag{4}$$

In the case $u_l > u_r$ one has $\Delta u<0$, therefore $\alpha>0$, as requested since the density $\rho$ cannot be $<0$. But in the case $u_l < u_r$, $\Delta u>0$, therefore $\alpha<0$, which is not acceptable for a density. Therefore the solution (4) is not physically acceptable in the case $u_l < u_r$ (one also finds it is unstable). Fortunately, in this case, one finds another solution, which is physically acceptable, [1,3,7]:

(5)

- if $x < u_l t$ then $u(x,t) = u_l$, $\rho(x,t) = \rho_l$ ( left-hand-side region),
- if $u_l t < x < u_r t$ then $u(x,t)$ undefined, $\rho(x,t) = 0$ (void region),
- if $x > u_r t$ then $u(x,t) = u_r$ , $\rho(x,t) = \rho_r$ (right-hand-side region).

This solution corresponds to the physics of the problem: in absence of pressure the two sides depart each other with their respective velocities.

## 3. Projection of delta waves

When a function is regular enough, say $L^\infty$, one usually projects it on a discretization lattice by taking its mean value on each cell [*ih-h/2, ih+h/2*]. This method lacks “continuity” when a delta peak is located close to an interface. Such a delta peak that, within the unavoidable uncertainty, would be located on the interface, could be as well attributed to any side. Numerical tests show that the trivial equal division to each side does not produce very good results. Then in presence of delta peaks, the knowledge of the function itself is not sufficient to permit a correct projection on a discretization lattice. The presence of delta peaks in the solution (3) of Riemann problems for the equations (S) of pressureless fluid dynamics and for the systems of physics in [6] therefore makes the projection step of a Godunov scheme non trivial. The delta peaks from the Riemann problems should have nontrivial right-hand-side and left-hand-side contributions to be discovered.

How can we treat the delta wave in the projection step of a Godunov scheme in the case $u_l > u_r$? The idea developed here is the following:

- In the case $u_l < u_r$, one applies the Godunov scheme using the solution (5), which has the usual form of step functions.
- Still in this case $u_l < u_r$, one seeks how to share the (non-physical) delta peaks in $\rho$ and $\rho u$ in (3), so as to obtain the same numerical scheme. The delta peak in $\rho$ is assumed to contribute to the left-hand-side cell by a factor $\lambda_l$ and to the right-hand-side cell by a factor $\lambda_r$, with $\lambda_l+\lambda_r=1$. Same for the delta peak in ($\rho u$), whose respective contributions are proportional to factors $\mu_l$ and $\mu_r$, with $\mu_l+\mu_r=1$. We compute explicitely the values $\lambda_l$, $\lambda_r$, $\mu_l$, $\mu_r$ that give the same numerical scheme as the one from the step functions solution.
- Now, in the case $u_l > u_r$, for each similar configuration of the waves, one adopts the same formulas for $\lambda_l$, $\lambda_r$, $\mu_l$, $\mu_r$ to share the delta peaks into left and right-hand-side contributions.

In the sequel of this section the sharing coefficients $\lambda_l, \lambda_r, \mu_l, \mu_r$ are calculated in the case $u_l < u_r$ as functions of the variables $u_l$, $u_r$, $\rho_l$, $\rho_r$. We denote by $\overline{w_l}$ (respectively $\overline{w_r}$) the mean value of a variable *w* on the segment [*-h/2,0*] (resp. [*0, h/2*]).

- Case $0 < u_l < u_r$. In this case $u_l < c < u_r$ from (4).

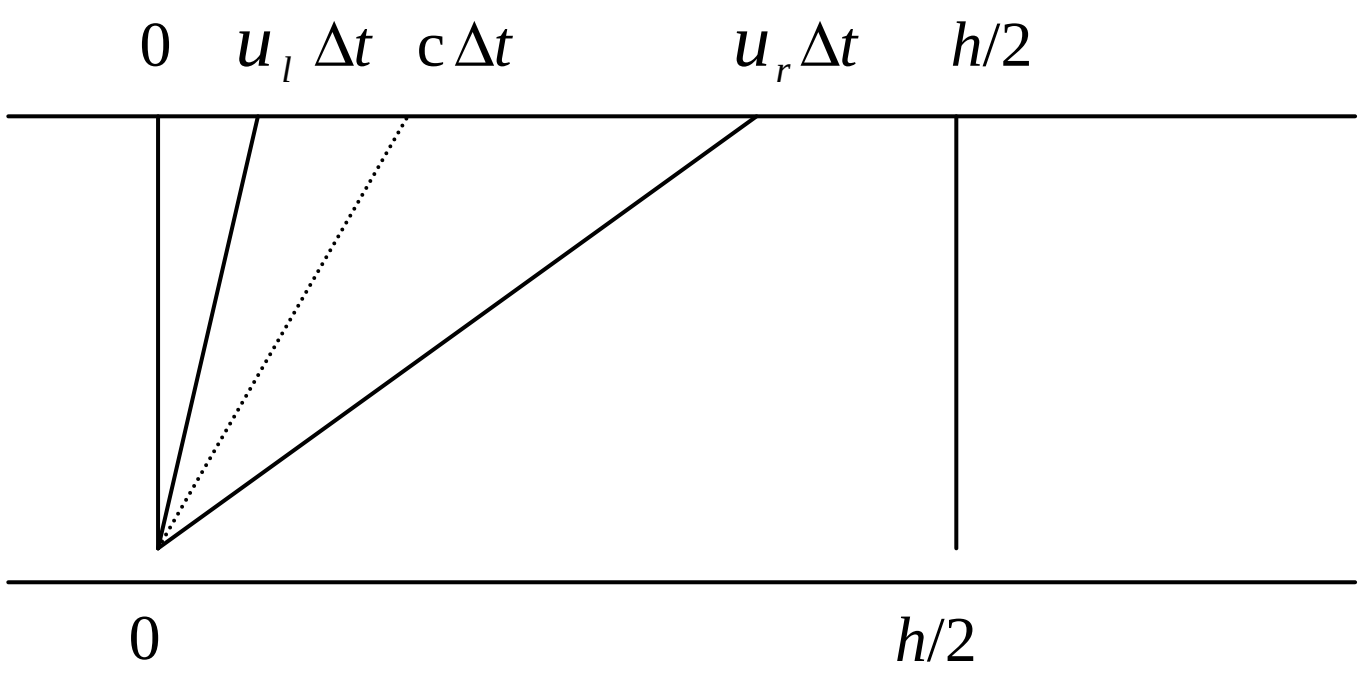

Projection of the step functions (two discontinuities of velocities $u_l$ and $u_r$, provided the CFL condition $u_r \Delta t < h/2$ i.e. $r\, u_r < 1/2$) gives from (5), as shown in the above picture:

$\overline{\rho_l} = \rho_l$,

$\overline{(\rho u)}_l = (\rho u)_l$,

$\overline{\rho_r} = (\rho_l u_l \Delta t + \rho_r (h/2 - u_r \Delta t))/(h/2) = \rho_r + 2r\rho_l u_l - 2r\rho_r u_r$,

$\overline{(\rho u)}_r = (\rho u)_r + 2\, r \rho_l\, u_l^{\,2} - 2r \rho_r\, u_r^{\,2}$.

Projection of the ( non-physical) delta wave gives from (3), (4):

$\overline{\rho_l} = (\rho_l h/2 + \lambda_l\, \alpha\, \Delta t)/(h/2) = \rho_l + 2\lambda_l\, \alpha\, r$,

$\overline{(\rho u)}_l = ((\rho u)_l h/2 + \mu_l\, \beta\, \Delta t)/(h/2) = (\rho u)_l + 2\mu_l\, \beta\, r$,

$\overline{\rho_r} = (\rho_l c \Delta t + \lambda_r\, \alpha\, \Delta t + \rho_r (h/2 - c\Delta t))/(h/2) = \rho_r + 2\rho_l cr + 2\lambda_r\, \alpha\, r - 2\rho_r cr$,

$\overline{(\rho u)}_r = ((\rho u)_l c\Delta t + \mu_r\, \beta\, \Delta t + (\rho u)_r (h/2 - c\Delta t))/(h/2) = (\rho u)_r + 2(\rho u)_l cr + 2\mu_r\, \beta\, r$
$-2(\rho u)_r\, cr$.

Identification of the two sets of formulas gives:

$\lambda_l = 0$,
$\mu_l = 0$,
$\rho_l u_l - \rho_r u_r = \rho_l c + \lambda_r\, \alpha - \rho_r c$,
$\rho_r u_r^{\,2} - \rho_l\, u_l^{\,2} = (\rho u)_l c + \mu_r\, \beta - (\rho u)_r\, c$.

Using (4), the last two formulas give, after immediate calculation, $\lambda_r = 1$ and $\mu_r = 1$. Therefore, in this case, the sharing coefficients are

(6) $\lambda_l = 0,\ \lambda_r = 1,\ \mu_l = 0,\ \mu_r = 1.$

This means that in this case the delta waves contribute only to the right-hand-side.

• Case $u_l < u_r < 0$. In this case one obtains similarly as above

(7) $\lambda_l = 1,\ \lambda_r = 0,\ \mu_l = 1,\ \mu_r = 0.$

• Case $u_l < 0 < u_r$ and $c > 0$.

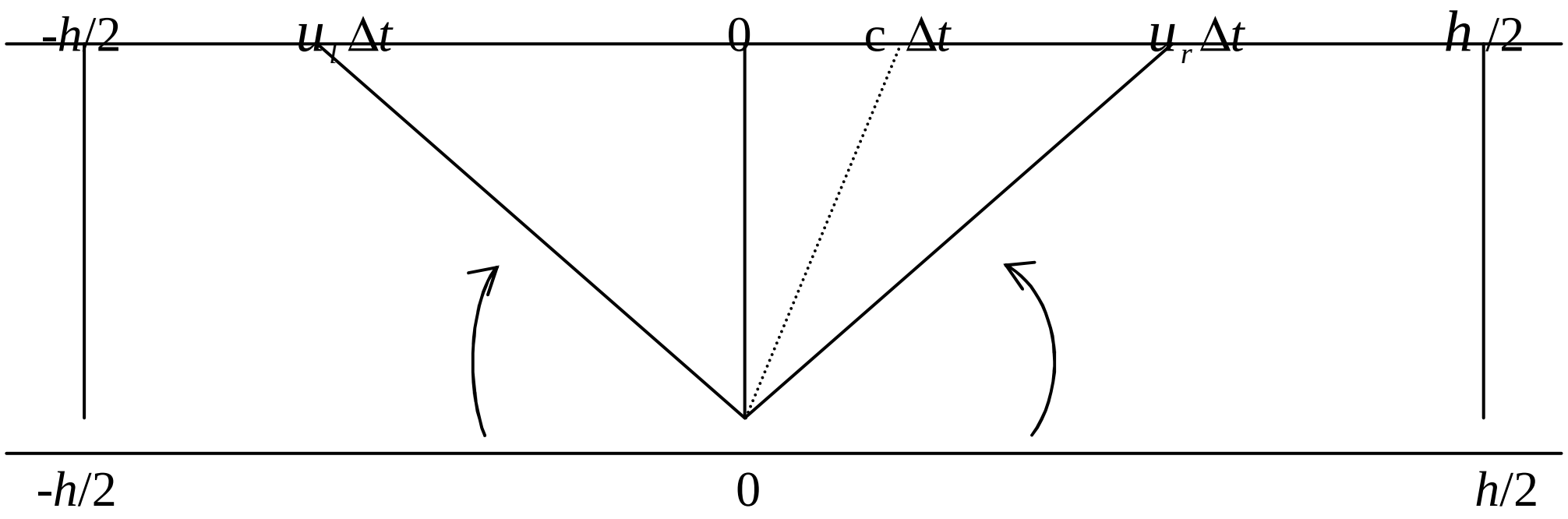

Projection of the step functions (provided the CFL condition $\max(|u_l|, u_r).\Delta t < h/2$ i.e. $r \max(|u_l|, u_r) < 1/2$) gives from (5):

$\overline{\rho_l} = \rho_l .(h/2 + u_l \Delta t)/(h/2) = \rho_l + 2\rho_l u_l r,$

$\overline{(\rho u)}_l = (\rho u)_l .(h/2 + u_l \Delta t)/(h/2) = (\rho u)_l + 2(\rho u)_l u_l r,$

$\overline{\rho_r} = \rho_r .(h/2 - u_r \Delta t)/(h/2) = \rho_r - 2\rho_r u_r r,$

$\overline{(\rho u)}_r = (\rho u)_r .(h/2 - u_r \Delta t)/(h/2) = (\rho u)_r - 2(\rho u)_r u_r r.$

Projection of the (non-physical) delta wave gives from (3), (4) exactly the same formulas as those obtained above in the case $0 < u_l < u_r$.

Identification of the two sets of formulas gives:

$\rho_l u_l = \lambda_l \alpha,$
$\rho_l u_l^2 = \mu_l \beta,$
$-\rho_r u_r = \rho_l c + \lambda_r \alpha - \rho_r c,$
$-(\rho u)_r u_r = (\rho u)_l c + \mu_r \beta - (\rho u)_r c.$

Thus one obtains the formulas for the left-hand-side and right-hand-side contributions of the delta wave:

(8) $\lambda_l = \rho_l u_l / \alpha,$

$\lambda_r = (-\rho_r u_r + c(\rho_r - \rho_l))/\alpha,$

$\mu_l = \rho_l u_l^2/\beta$,
$\mu_r = (-\rho_r u_r^2 + c(\rho_r u_r - \rho_l u_l))/\beta$.

- Case $u_l < 0 < u_r$ and $c<0$. Similar calculations give:

(9) $\lambda_l = (\rho_l u_l + c(\rho_r - \rho_l))/\alpha$,
$\lambda_r = -\rho_r u_r/\alpha$,
$\mu_l = (\rho_l u_l^2 + c(\rho_r u_r - \rho_l u_l))/\beta$,
$\mu_r = -\rho_r u_r^2/\beta$.

Sum up of these results: rule of splitting of the delta peak observed in the case $u_l < u_r$. The splitting of the delta peaks into a left-hand-side contribution and a right-hand-side contribution depends on the (left or right hand-side) positions of the three waves under concern: those of velocities $u_l$, $u_r$ and the delta peak of velocity $c$. Looking at the above four cases in which $u_l < u_r$, one arrives at the conclusion that the following rule always hold to evaluate the $\lambda_l$, $\lambda_r$ factors in the contribution of the delta peak in $\rho$:

- $\lambda$-contribution to the side where the wave of velocity $u_r$ is located: $-\rho_r u_r/\alpha$,
- $\lambda$-contribution to the side where the wave of velocity $u_l$ is located : $\rho_l u_l/\alpha$,
- $\lambda$-contribution to the side where the delta peak is located: $c(\rho_r - \rho_l)/\alpha$.

Note that the contributions are null if $u_r=0$, $u_l=0$ and $c=0$: there is no ambiguity when a wave lies at the interface.

Example: in the case $u_l < 0 < u_r$ and $c>0$,
the wave of velocity $u_r$ contributes to the right-hand-side, i.e. to $\lambda_r$, by $-\rho_r u_r/\alpha$,
the wave of velocity $c$ contributes to the right-hand-side, i.e. to $\lambda_r$, by $c(\rho_r - \rho_l)/\alpha$,
the wave of velocity $u_l$ contributes to the left-hand-side, i.e. to $\lambda_l$, by $\rho_l u_l/\alpha$.
Sum up of all contributions: $\lambda_r = -\rho_r u_r/\alpha + c(\rho_r - \rho_l)/\alpha$ and $\lambda_l = \rho_l u_l/\alpha$. One recovers (8).

For the $\mu_l, \mu_r$ factors in the contribution of the delta peak in $(\rho u)$ the rule is:

- $\mu$-contribution to the side where the wave of velocity $u_r$ is located: $-\rho_r u_r^2/\beta$,
- $\mu$-contribution to the side where the wave of velocity $u_l$ is located : $\rho_l u_l^2/\beta$,
- $\mu$-contribution to the side where the delta peak is located: $c(\rho_r u_r - \rho_l u_l)/\beta$.

Then, the method of this paper consists in adopting this rule (obtained in the known case $u_l < u_r$) in the (unknown) case $u_r < u_l$ for the splitting of delta peaks. How can it be justified? One could think that the proper formulas for the projection of delta waves are the same whether they are physical or non-physical. An *a posteriori* justification is provided by the convergence result in theorem 3, although it does not cover all main cases, and the convergence result in theorem 4 below, which concerns a very slight modification of the scheme but covers all cases. To validate the scheme and test its quality one has compared the numerical results with other solutions (figures 1, 2 and 3): one obtains the correct solution in all cases.

Remark. This method can be applied to nonlinear systems, which are systems of physics or systems obtained from a splitting of systems of physics and that are variants of the system of pressureless fluid dynamics. Numerical approximations are proposed in [6] in the Newtonian and relativistic cases, where their use permits to observe the main steps of large structure formation according to the cosmic epochs.

## 4. Interpretation of the splitting rule

In the case $u_r > u_l$, in which the formulas were obtained, one has a void region separated by discontinuities of velocities $u_l$ and $u_r$. In the case $u_r < u_l$ one has instead some phenomenon looking intuitively like a collision of two volumes of fluid. Does the splitting rule adopted allow an intuitive interpretation in the collision case?

**Theorem 1:** In the "collision case" the splitting rule can be interpreted as an "interpenetration" of the two volumes of fluid. It amounts to decompose the physical phenomenon that occurs at the level of "infinitesimal cells", into first "free streaming through the interfaces of cells limited to small volume by the CFL condition", then "mixing" inside each cell (the mean value) so as to have well defined density and momentum in each cell.

Proof of the theorem:

- Case $0 < u_r < u_l$.

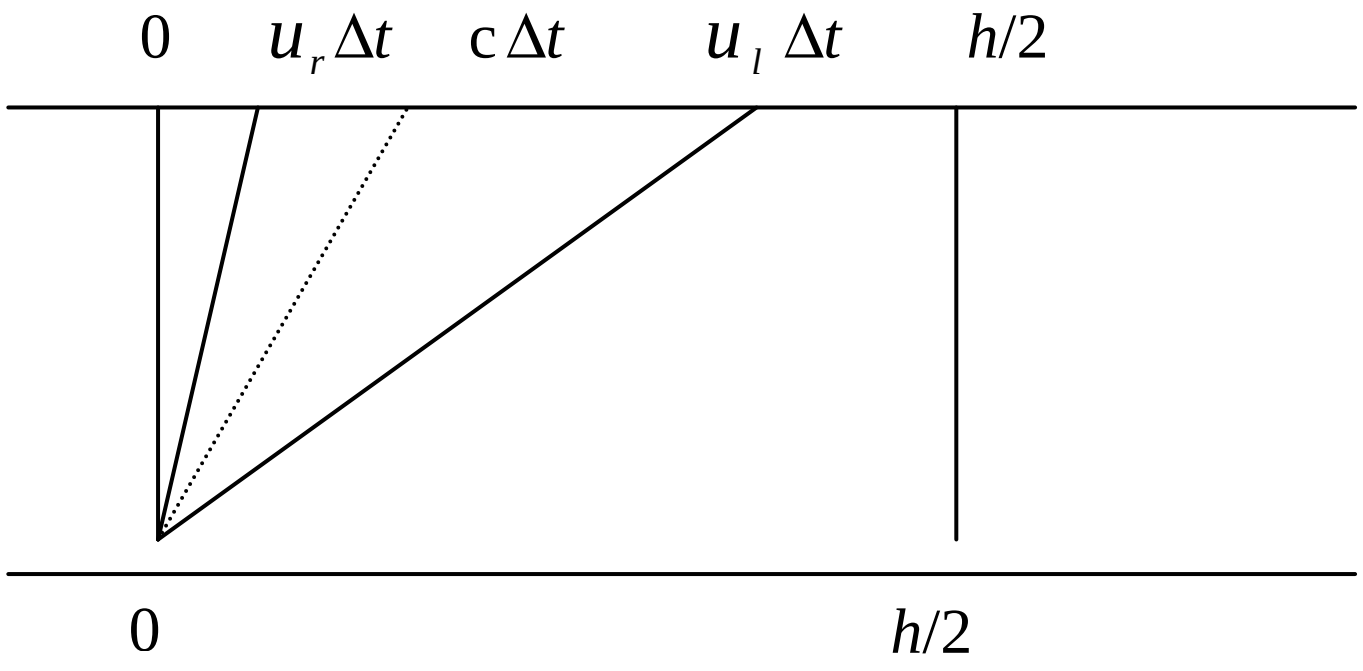

From (6) $\overline{\rho_l} = \rho_l$ and $\overline{\rho_r} = \rho_r + 2r\rho_l u_l - 2r\rho_r u_r$. The projection in case of interpenetration of the two fluids gives $\overline{\rho_l} = \rho_l$ and $\overline{\rho_r} = (u_l \Delta t\, \rho_l + (h/2-u_r\Delta t)\,\rho_r)/(h/2) = \rho_r + 2r\rho_l u_l - 2r\rho_r u_r$. Same formulas hold for $\overline{(\rho u)}_l$ and $\overline{(\rho u)}_r$. Thus one recovers the formulas from (6). Same results are obtained in the case $u_r < u_l < 0$.

- Case $u_r < 0 < u_l$ and $c > 0$.

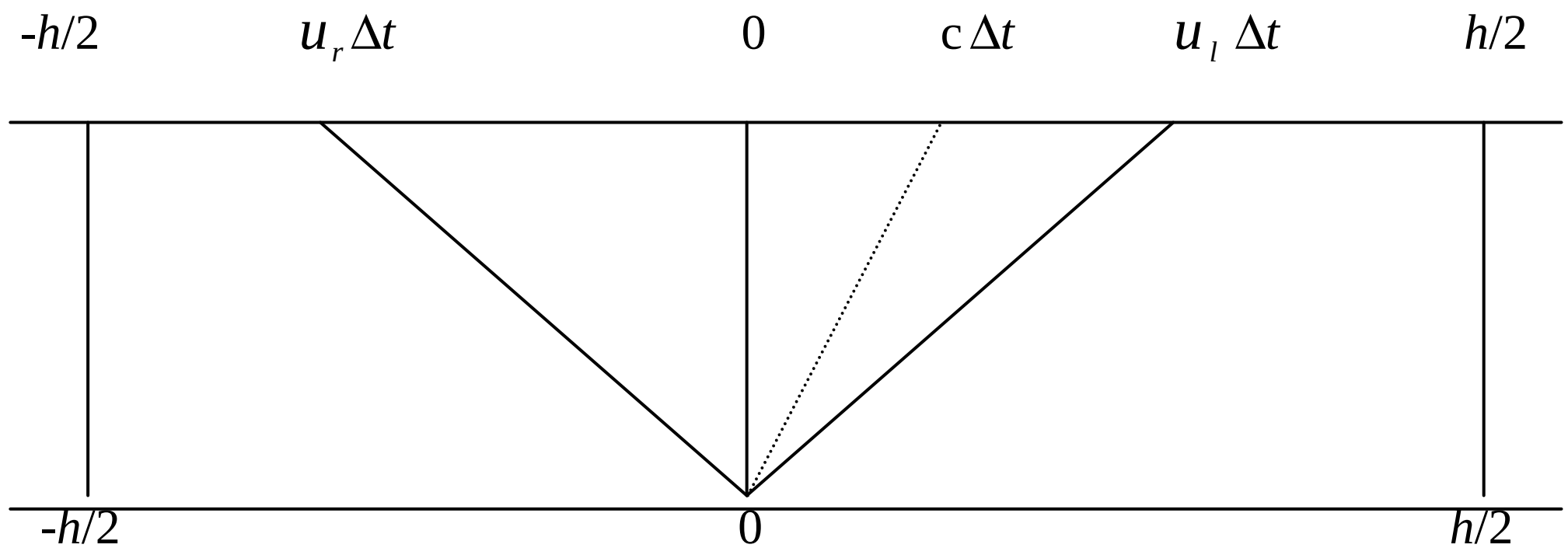

In this case the splitting rule gives $\lambda_l = -\rho_r u_r/\alpha$, $\lambda_r = (\rho_l u_l + c(\rho_r - \rho_l))/\alpha$, $\mu_l = -\rho_r u_r^2/\beta$, $\mu_r = (\rho_l u_l^2 + c(\rho_r u_r - \rho_l u_l))/\beta$.

Then, from the splitting rule

$\overline{\rho_l} = ((h/2)\,\rho_l + \lambda_l\, \alpha\, \Delta t)/(h/2) = \rho_l + 2r\lambda_l\, \alpha = \rho_l + 2r(-\rho_r u_r)$,

$\overline{\rho_r} = ((h/2) - c\Delta t)\,\rho_r + c\Delta t\, \rho_l + \lambda_r\, \alpha\, \Delta t)/(h/2) = \rho_r + 2rc(\rho_l - \rho_r) + 2r(\rho_l u_l + c(\rho_r - \rho_l)) = \rho_r + 2r\rho_l u_l$.

$\overline{(\rho u)}_l = ((h/2)\,(\rho u)_l + \mu_l\, \beta\, \Delta t)/(h/2) = (\rho u)_l + 2r\mu_l\, \beta = (\rho u)_l + 2r(-\rho_r u_r^2)$,

$\overline{(\rho u)}_r = (((h/2) - c\Delta t)(\rho u)_r + c\Delta t(\rho u)_l + \mu_r\, \beta\, \Delta t)/(h/2) = (\rho u)_r + c\; 2r((\rho u)_l - (\rho u)_r) + 2r(\rho_l u_l^2 + c(\rho_r u_r - \rho_l u_l)) = (\rho u)_r + 2r\rho_l u_l^2$.

The projection in case of interpenetration gives:

$\overline{\rho_l} = (-u_r\Delta t\, \rho_r + (h/2)\, \rho_l)/(h/2) = \rho_l - 2r\rho_r u_r$,

$\overline{\rho_r} = ((h/2)\, \rho_r + u_l\, \Delta t\, \rho_l)/(h/2) = \rho_r + 2\rho_l u_l$.

$\overline{(\rho u)}_l = (-u_r\Delta t(\rho u)_r + (h/2)(\rho u)_l)/(h/2) = (\rho u)_l - 2r\rho_r u_r^2$,

$\overline{(\rho u)}_r = ( (h/2)\ \rho_r\ u_r + u_l\ \Delta t\ \rho_l u_l)/(h/2) = \rho_r\ u_r + 2r\rho_l u_l^{\ 2}$ .

We proceed the same way in the case $u_r<0<u_l$ and $c<0$.

Remark 2. In [6] more complicated systems of physics are presented, for which physical intuition cannot give the numerical scheme, and that can be treated in the same way.

Remark 3. The scheme can be written:

$$\overline{\rho_l} = [\,\rho_l \int_0^{h/2} [-h/2, rhu_l](x)dx \ + \rho_r \int_0^{h/2} [rhu_r, h/2](x)dx\,]/(h/2),$$

$$\overline{\rho_r} = [\,\rho_l \int_{-h/2}^{0} [-h/2, rhu_l](x)dx \ + \rho_r \int_{-h/2}^{0} [rhu_r, h/2](x)dx\,]/(h/2),$$

and similar formulas for $\overline{(\rho u)}_l$ and $\overline{(\rho u)}_r$, with the notation $[a,b](x)=1$ if $a<x<b$, $0$ if not. In case $u_l< u_r$ this follows from the classical Godunov scheme. The splitting rule consists in adopting the same formulas in the unknown case $u_l> u_r$ .

## 5. Comparison with the scheme in [1]

**Proposition1:** If $u_{i-1}^n>0$, $u_i^n>0$, $u_{i+1}^n>0$ then the scheme in [1] and the scheme in this paper give the same values of $\rho_i^{n+1}$ and $u_i^{n+1}$. The same result holds in case $u_{i-1}^n<0$, $u_i^n<0$, $u_{i+1}^n<0$.

Proof: From (1) $u_{i-1/2}^n = u_{i-1}^n$, $\rho_{i-1/2}^n = \rho_{i-1}^n$, $u_{i+1/2}^n = u_i^n$, $\rho_{i+1/2}^n = \rho_i^n$. Therefore $\rho_i^{n+1} = \rho_i^n - r\rho_i^n\ u_i^n + r\rho_{i-1}^n\ u_{i-1}^n$, $\rho u_i^{n+1} = \rho u_i^n - r\rho u_i^n u_i^n + r\rho u_{i-1}^n u_{i-1}^n$. Now, from theorem 1, one checks easily that

$\rho_i^{n+1} = (u_{i-1}^n\ \Delta t\ \rho_{i-1}^n + (h - u_i^n\ \Delta t)\rho_i^n)/h = r\rho_{i-1}^n\ u_{i-1}^n + (1 - r\,u_i^n)\rho_i^n$,

which is the same formula. The calculations give the same result in case of $\rho u$.

Looking at the various cases after (1) the scheme in this paper differs from the scheme in [1] in the case $u_i^n>0$ and $u_{i+1}^n<0$ which is a “collision case”.

## 6. Numerical schemes in 1D, 2D and 3D

Let $a$ and $b$ be two numbers with $a<b$, $|a|<1/2$, $|b-1|<1/2$. Set

$$L(a,b) = \text{length of } [0,1] \cap [a,b] = \min(1,b) - \max(0,a).$$

In both cases (void region if $u_r > u_l$ and interpenetration if $u_r < u_l$) the volume of matter which lies in [*ih-h/2, ih+h/2*] at time $t_n = n.\Delta t = nrh$ is transported between the times $t_n$ and $t_{n+1}$ with velocity $u_i^n$. Then after transport, from the CFL condition below, this volume is to be found at time $t_{n+1}$ in the three cells [*ih-3h/2, ih-h/2*], [*ih-h/2, ih+h/2*] and [*ih+h/2, ih+3h/2*] as well as in the segment [$ih-h/2+ u_i^n \Delta t$, $ih+h/2+ u_i^n \Delta t$]. Therefore the mean values $\rho_i^{n+1}$ and $(\rho u)_i^{n+1}$ in the cell [*ih-h/2, ih+h/2*] at time $t_{n+1}$ are given by the following formulas that sum up the possible contributions from the cells *i-1, i* and *i+1*:

$$\rho_i^{n+1} = \rho_{i-1}^n L(-1+ru_{i-1}^n, ru_{i-1}^n) + \rho_i^n L(ru_i^n, 1+ru_i^n) + \rho_{i+1}^n L(1+ru_{i+1}^n, 2+ru_{i+1}^n),$$

(10)

$$(\rho u)_i^{n+1} = (\rho u)_{i-1}^n L(-1+ru_{i-1}^n, ru_{i-1}^n) + (\rho u)_i^n L(ru_i^n, 1+ru_i^n) + (\rho u)_{i+1}^n L(1+ru_{i+1}^n, 2+ru_{i+1}^n),$$

which are completed by

$$u_i^{n+1} = \frac{(\rho u)_i^{n+1}}{\rho_i^{n+1}}.$$

The CFL condition is:

$$r \max\{|u_i^n|\} \leq 1.$$

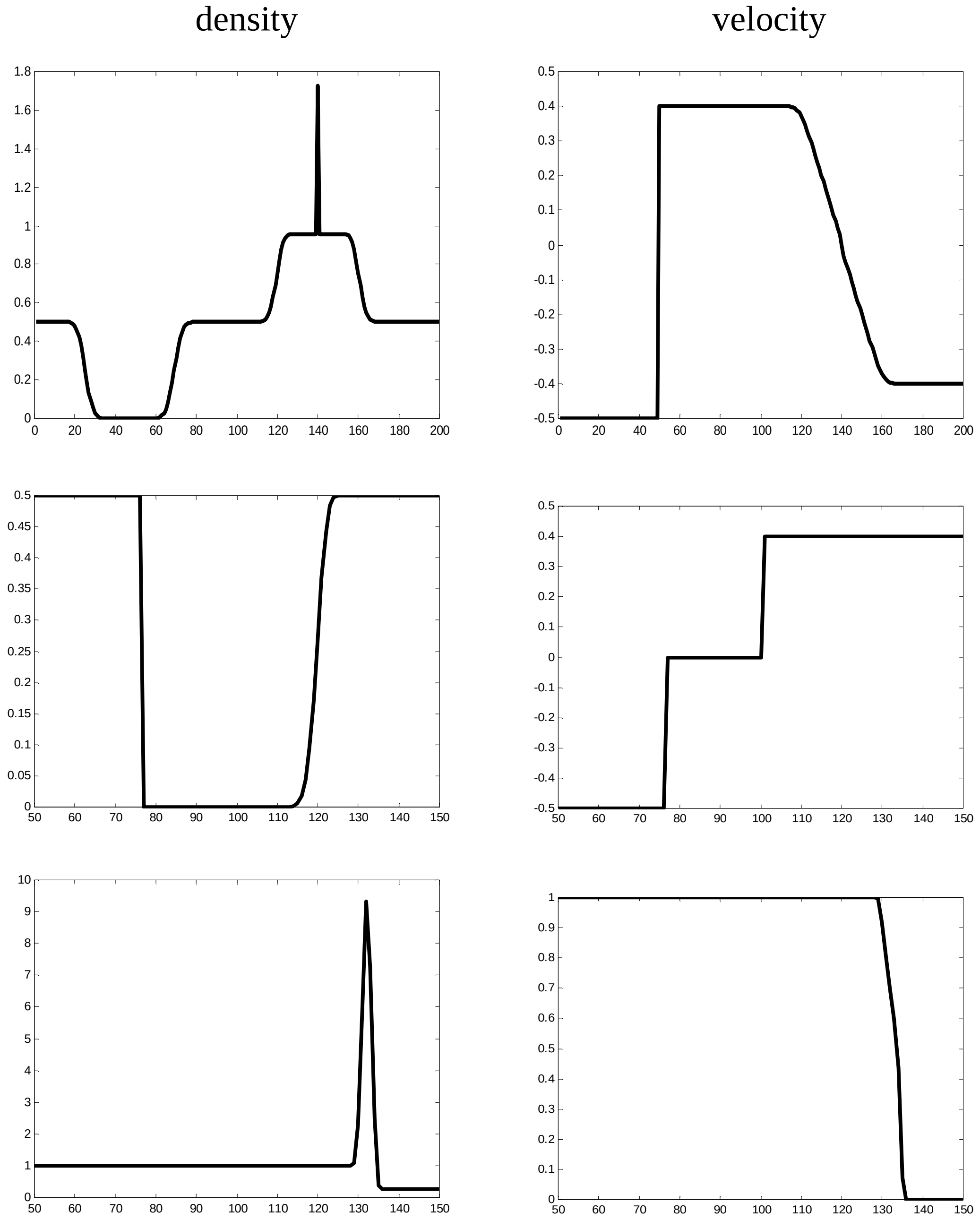

**figure 1: from top to bottom the results of the three 1D tests in [2] using the scheme presented here. One obtains results practically identical to those in [2] (200 points mesh, same time) and to the exact solution given in [2] (note that the velocity is meaningless in the void regions). The minor defect (the spike) observed in top-left is identically observed in the schemes tested in [2]. In the second test the minimum value of the vacuum state is about $10^{-300}$. Concerning the third test it has been observed that the presence of gravitation in the context of [6] reduces considerably the support of delta peaks.**

Collision of two dust clouds

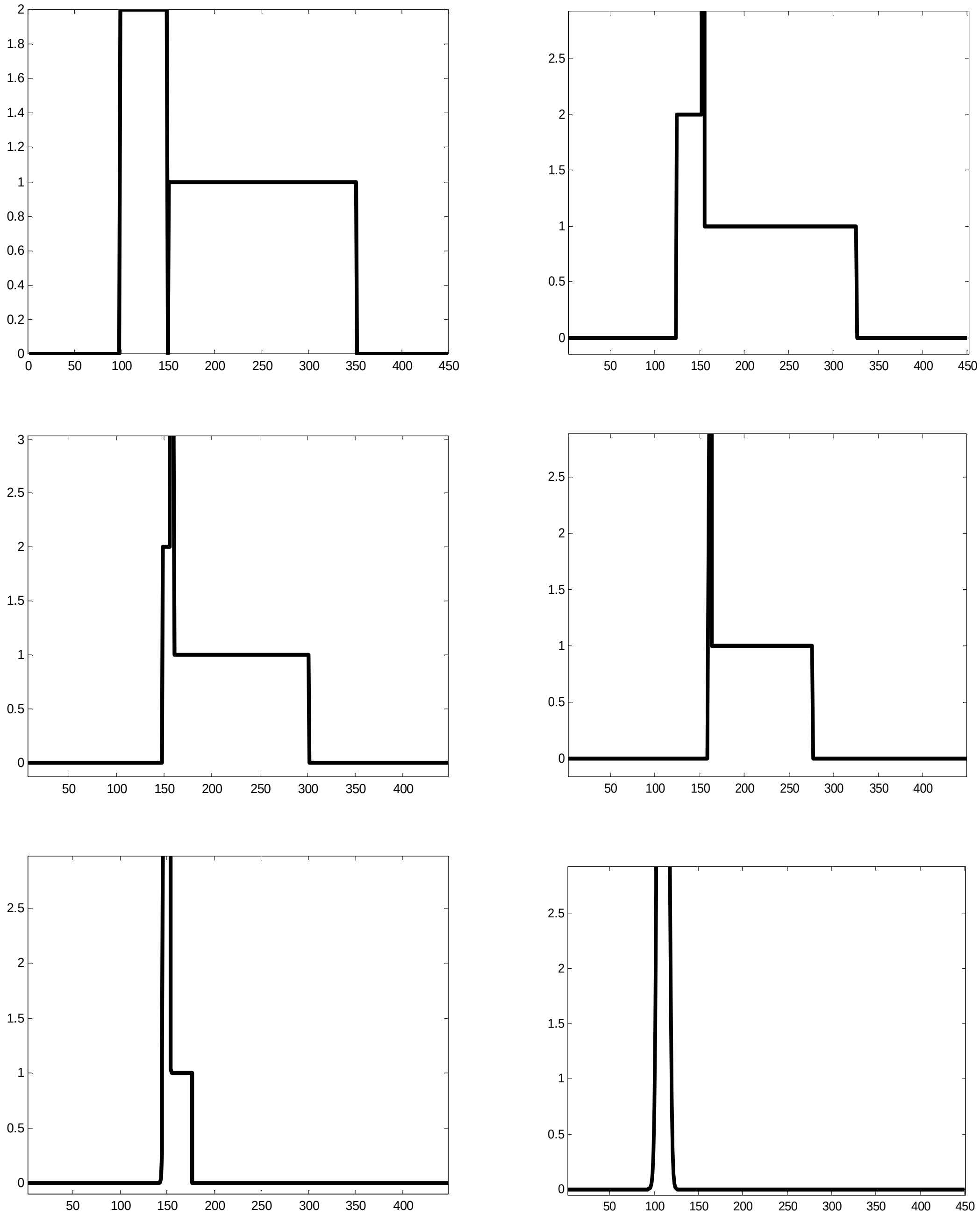

**figure 2: A numerical test from [7]. Two clouds of pressureless gas collide at t=0. From left to right and top to bottom t = 0, 0.5, 1, 1.5, 3.5, 6. The space has been divided into 450 meshes and r=1. One observes complete agreement with the numerical results and the exact solution given in [7]. At time t=6 all mass is concentrated in a single delta peak which is smeared out over 24 cells (this defect that "double delta peaks" adjacent to vacuum states are not captured as sharply as expected has been observed also in [7]). Gravitation is present in the context of cosmology [6] and it will repair this defect since it has been observed it can reduce the support of delta peaks to one cell even in 2D and 3D tests. The minimum value of the vacuum states is less than $10^{-300}$. The computation time on a standard PC has been 0.24 second.**

numerical solutions comparison exact/numerical

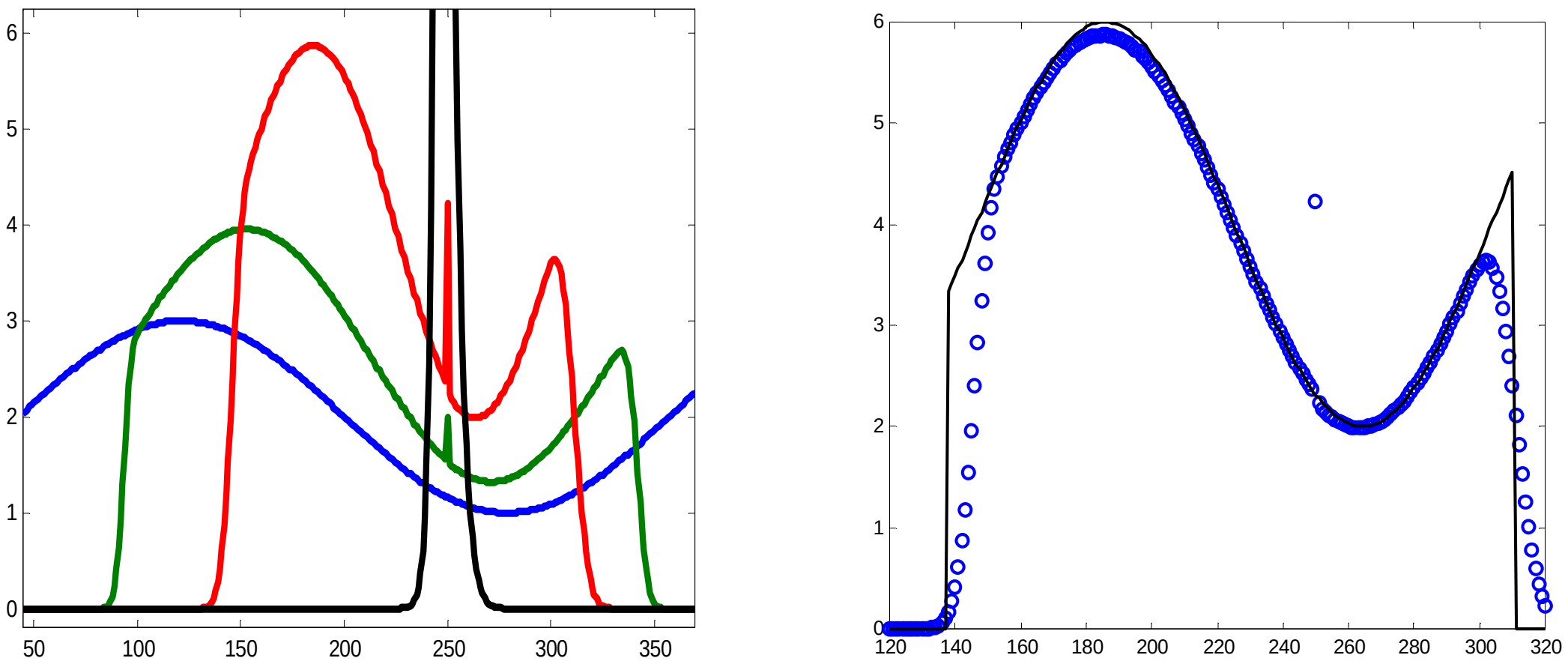

**Figure 3. The 1D test 4 in [4]: the velocity changes its sign in a region with varying density (according to [4, top of p3] many methods fail in this test). Left figure: blue at t=0, green at t=0.25, red at t=0.5, black =delta peak solution at t=1. The final collapse of the numerical solution occurs exactly at its correct time and location (t=1, x=1 in units of [4]) in form of a delta peak whose base encompasses 35 cells. Right figure: at t = 0.5 exact solution (black) and numerical solution (blue 'o'); one observes numerical dissipation on the discontinuities and an isolated point. For t>1 the delta peak moves to the right as expected but its base enlarges; gravitation which is present in the applications in [6] will repair this defect. The computation time on a standard PC has been 0.1 second.**

In the 2D case the equations of pressureless fluid dynamics are

$$\rho_t + (\rho u)_x + (\rho v)_y = 0,$$
$$(11) \qquad (\rho u)_t + (\rho u^2)_x + (\rho uv)_y = 0,$$
$$(\rho v)_t + (\rho uv)_x + (\rho v^2)_y = 0,$$

in which $(u,v)$ is the velocity vector. The 2D space $(x,y)$ is divided into square cells of side $h$ and center $(ih, jh)$, $i,j \in \mathbb{Z}$,

$$C_{i,j} = \{(x,y) \text{ such that } ih-h/2 < x < ih+h/2 \text{ and } jh-h/2 < y < jh+h/2\}.$$

One considers the eight cells that have a joint boundary with the cell $C_{i,j}$ . The numerical scheme is the passage from the set $\left\{\rho^n_{i,j}, (\rho u)^n_{i,j}, (\rho v)^n_{i,j}\right\}_{i,j \in Z}$ to the set $\left\{\rho^{n+1}_{i,j}, (\rho u)^{n+1}_{i,j}, (\rho v)^{n+1}_{i,j}\right\}_{i,j \in Z}$. Analogously to the 1D case ( it is clear that system

(11) amounts to a translation by the vector ($u.\Delta t$, $v.\Delta t$) in place of $u.\Delta t$, as found in the 1D case), set:

*A( a, b)* = area of the intersection of the square of vertices *{(0,0), (0,1), (1,0), (1,1)}* with the square of vertices *{(a, 1+a), (a, b), (b, 1+b), (1+a,1+b)}*. One has:

$$A(a, b) = L(a,1+a).\ L(b,1+b).$$

The square $C_{p,q}$ is translated by the vector ($u^n_{p,q}\,\Delta t$, $v^n_{p,q}\,\Delta t$). Taking into account the cell $C_{i,j}$ itself and its eight neighbours, one obtains the formulas of the 2D scheme, if $\omega = \rho$, $\rho$u, and $\rho$v successively:

$$(12) \qquad \omega^{n+1}_{i,j} = \sum_{\substack{\lambda=i-1,i,i+1\\ \mu=j-1,j,j+1}} \omega^n_{\lambda,\mu}\, A(\lambda - i + ru^n_{\lambda,\mu}, \mu - j + rv^n_{\lambda,\mu}),$$

completed by

$$u^{n+1}_{i,j} = \frac{(\rho u)^{n+1}_{i,j}}{\rho^{n+1}_{i,j}}, \quad v^{n+1}_{i,j} = \frac{(\rho v)^{n+1}_{i,j}}{\rho^{n+1}_{i,j}}.$$

The CFL condition is

$$r.\max\{|u^n_{i,j}|, |v^n_{i,j}|\} \le 1.$$

Adaptation of the scheme to 3D is immediate: it suffices to define an auxiliary function *V* (=”volume”) similar to *A* (=”area”), and an auxiliary function taking into account the 26 neighbours plus the cell itself, see a 3D numerical scheme in [6]. The 2D and 3D tests have worked quite well, like the 1D ones given above. The numerical results approximate the exact solution when it exists; some numerical viscosity can be observed on discontinuities and delta shocks, like in all schemes of the same nature (but not in the sticky particule method of [4]).

## 7. Stability of the scheme and weak*compactness

**Theorem 2:** The scheme is $L^1$ stable in the variables $\rho$ and $\rho u$, and $L^\infty$ stable in the variable *u* for which the maximum principle holds. These results hold in any space dimension.

Proof: From (10) (and eliminating boundary effects, for instance if $\rho$ and $\rho u$ are null close to the left-hand-side and right-hand-side boundaries) one has:

$$\sum_{i\in Z}\rho_i^{n+1} = \sum_{i\in Z}\rho_i^n = \sum_{i\in Z}\rho_i^0$$

and same for the variable $\rho$u. Since $\rho$ is positive this implies the $L^1$ stability of the scheme in the $\rho$ variable.
Now assume that

$$a \le u_i^n \le b \ \forall i \in Z .$$

Since $\rho_i^n \ge 0$ this implies that

$$a\,\rho_i^n \le \rho_i^n u_i^n \le b\,\rho_i^n \ \ \forall i \in Z .$$

Since the quantities L in (10) are positive, the second line in (10) implies that

$$a\,\rho_i^{n+1} \le (\rho\, u)_i^{n+1} \le b\,\rho_i^{n+1} \ \ \forall i \in Z .$$

Assuming $\rho_i^{n+1} > 0$ (if not $u_i^{n+1}$ is meaningless) this implies that

$$a \le u_i^{n+1} \le b \ \forall i \in Z .$$

which is the maximum principle in the variable $u$. Now

$$\sum_{i\in Z} |\rho\, u|_i^n \ \le \ \max(|a|,|b|)\sum_{i\in Z}\rho_i^n$$

implies that the scheme is $L^1$ stable in the $\rho u$ variable. This permits to use compactness in the usual way. The scheme takes place in a finite interval I. Embedding $L^1$(I) into the Banach space M of Radon measures on I, there exists a sequence $h_p \to 0$ such that the corresponding sequences $(\rho_p)$ and $(u_p)$ converge in the weak-* topologies in M and $L^\infty$ respectively to a measure $\rho$ and to a $L^\infty$ function $u$. The proof can be reproduced without modifications, other than the obvious ones, in any space dimension.

## 8. A convergence result in case the velocity has constant sign

**Theorem 3**: In the case velocity has a constant sign the scheme converges to a pair ($\rho$,$u$) of a measure $\rho$ and a function $u$ solution in the sense of distributions of the equations of pressureless fluid dynamics on $IR^n \times [0,+\infty)$ when the initial condition is $L^1$ in $\rho$ and $L^\infty$ in $u$. Note that this holds in any space dimension.

**Corollary:** In the case the velocity has a constant sign the 1D-scheme in [1] is convergent (from proposition 1).

Proof of theorem 3: Set $IR^{2+}=\{x\in IR,\ t>0\}$ and $D(IR^{2+})$ = the space of all $C^\infty$ functions $IR^{2+}\mapsto IR$ with compact support in $IR^{2+}$. We are going to prove that $\forall\ \psi\in D(IR^{2+})$

$\int \{\rho\psi_t + (\rho u)\psi_x\}\,dxdt$ and $\int \{(\rho u)\psi_t + (\rho u^2)\psi_x\}\,dxdt$

tend to 0 as $h\to 0$, if $\rho$ and $u$ are the step functions from the numerical scheme with space meshsize $h$ and time meshsize $rh$ (fixed $r$). To this end we do as follows. From the discretization

$$\int \rho\psi_t dxdt = \sum_{i,n} \rho_i^n \Big(\int_{(i-1/2)h<x<(i+1/2)h,\,nrh<t<(n+1)rh} \psi_t(x,t)dxdt\Big).$$

$\psi_t(x,t) = \psi_t(ih,nrh) + rh\,\psi_{tt}(\ ) + h\psi_{tx}(\ ) = (\psi_t)_i^n + h\,O(1)$ where ( ) concerns intermediate points and $(\psi_t)_i^n = \psi_t(ih,nrh)$. Then

$$\int \rho\psi_t dxdt = \sum_{i,n} \rho_i^n rh^2 (\psi_t)_i^n + \sum_{i,n} \rho_i^n\, rh^2 h\,O(1) = rh^2 \sum_{i,n} \rho_i^n (\psi_t)_i^n + h\,O(1)$$

from the $L^1$ stability in $\rho$ and since $\sum_n rh \le$ fixed value $T$ due to the compact support of $\psi$. A similar bound holds for $\int \rho u\,\psi_x\,dxdt$ from the $\rho u$ $L^1$ stability.

Therefore, setting $(\psi_x)_i^n=\psi_x(ih,nrh)$,

(13) $\int \{\rho\psi_t + (\rho u)\psi_x\}\,dxdt = rh^2 \sum_{i,n} \{\rho_i^n (\psi_t)_i^n + (\rho u)_i^n\, (\psi_x)_i^n\} + h\,O(1).$

Similarly

(14) $\int \{(\rho u)\psi_t + (\rho u^2)\psi_x\}\,dxdt =$

$$rh^2 \sum_{i,n} \{(\rho u)_i^n (\psi_t)_i^n + (\rho u^2)_i^n\, (\psi_x)_i^n\} + h\,O(1).$$

Now the formulas $(\psi_t)_i^n = (\psi_i^{n+1}-\psi_i^n)/rh + h\,O(1)$, $(\psi_x)_i^n=(\psi_{i+1}^n-\psi_i^n)/h + h\,O(1)$ permit to rewrite (13) as

$$\int \{\rho\psi_t + (\rho u)\psi_x\}\,dxdt = rh^2 \sum_{i,n} \{\rho_i^n (\psi_i^{n+1}-\psi_i^n)/rh + (\rho u)_i^n (\psi_{i+1}^n-\psi_i^n)/h\}$$

$$+ rh^2 \sum_{i,n} \rho_i^n hO(1) + rh^2 \sum_{i,n} (\rho u)_i^n hO(1) + hO(1).$$

From the $L^1$ stability in $\rho$, $\rho u$ and since $\sum_n rh \le T$):

$$hO(1) \sum_{i,n} rh^2 \rho_i^n + hO(1) \sum_{i,n} rh^2 |\rho u|_i^n \le hO(1).$$

Therefore

(15) $\int \{\rho\psi_t + (\rho u)\psi_x\}\,dxdt = h\sum_{i,n} \{\rho_i^n (\psi_i^{n+1}-\psi_i^n) + r(\rho u)_i^n (\psi_{i+1}^n-\psi_i^n)\} +$

$hO(1) = h\{\sum_{i,n} \rho_i^{n-1}\psi_i^n - \sum_{i,n} \rho_i^n\psi_i^n + r\sum_{i,n} (\rho u)_{i-1}^n\, \psi_i^n - r\sum_{i,n} (\rho u)_i^n \psi_i^n\} + hO(1)$

$$= -h\sum_{i,n} \{\rho_i^n - \rho_i^{n-1} + r[(\rho u)_i^n - (\rho u)_{i-1}^n]\}\, \psi_i^n + hO(1).$$

After another change of indices

$$\sum_{i,n} [(\rho u)_i^n - (\rho u)_{i-1}^n]\} \ \psi_i^n = \sum_{i,n} \ [(\rho u)_i^{n-1} - (\rho u)_{i-1}^{n-1}]\} \ (\psi_i^n + \psi_i^{n-1} - \psi_i^n).$$

Claim: $\sum_{i,n} \ [(\rho u)_i^{n-1} - (\rho u)_{i-1}^{n-1}]\} \ (\psi_i^{n-1} - \psi_i^n) = h\mathrm{O}(1).$

Proof of the claim: the left-hand-side member is equal to $\sum_{i,n} \ (\rho u)_i^{n-1}\{(\psi_i^{n-1} - \psi_i^n) + (\psi_{i+1}^n - \psi_{i+1}^{n-1})\}$. The factor { } $= h^2\,\mathrm{O}(1)$. Then use the $L^\infty$ and $L^1$ stability in $u$ and $\rho$ respectively.

End of the proof of theorem 3:

Therefore $\sum_{i,n} [(\rho u)_i^n - (\rho u)_{i-1}^n]\}\psi_i^n = \sum_{i,n} \ [(\rho \mathrm{u})_i^{n-1} - (\rho u)_{i-1}^{n-1}]\}\psi_i^n + h\mathrm{O}(1).$

Using (15)

(16) $\int \{\rho\, \psi_t + (\rho u)\psi_x\}\mathrm{d}x\mathrm{d}t =$

$$-h\sum_{i,n} \{\rho_i^n - \rho_i^{n-1} + r[(\rho u)_i^{n-1} - (\rho u)_{i-1}^{n-1}]\}\psi_i^n + h\mathrm{O}(1).$$

Similarly

(17) $\int \{(\rho u)\psi_t + (\rho u^2)\psi_x\}\mathrm{d}x\mathrm{d}t =$

$$-h\sum_{i,n} \{(\rho u)_i^n - (\rho u)_i^{n-1} + r[(\rho u^2)_i^{n-1} - (\rho u^2)_{i-1}^{n-1}]\}\psi_i^n + h\mathrm{O}(1).$$

Now let us assume that $u_i^n \geq 0 \ \forall\, i$ and $\forall\, n$. If $w = \rho$ or $\rho u$ then one has from the scheme

$$w_i^n = ru_{i-1}^{n-1}w_{i-1}^{n-1} + (1-ru_i^{n-1})w_i^{n-1} = w_i^{n-1} - r(w_i^{n-1}u_i^{n-1} - w_{i-1}^{n-1}u_{i-1}^{n-1}).$$

Inserting this formula with $w=\rho$ into (16) and $w = \rho u$ into (17) the first terms in the second members disappear. The same proof holds if all velocities are $\leq 0$. Clearly the same proof holds in 2D and 3D (after (17) the scheme gives various terms: those of smaller order are only relevant and correspond to the terms in (16)(17)).

Question: Is the scheme still convergent without this strong assumption on the velocity? In the next section the convergence of a very slightly modified scheme is proved. The slight modification brings a simplification of the proof of convergence: it amounts to solve a system which has same solution modulo a translation, and for which the same proof as the one of theorem 3 applies because all velocities are positive.

## §9. A convergence result for a slightly modified scheme

Let us consider a system of conservation laws:

(18) $$\frac{\partial}{\partial t}u(x,t)+\frac{\partial}{\partial x}f(u(x,t))=0.$$

For $c > 0$ set

(19) $$U(x,t)=u(x-ct,t).$$

Then

$$\frac{\partial}{\partial t}U(x,t)=-cD_1u(x-ct,t)+D_2u(x-ct,t),$$
$$\frac{\partial}{\partial x}f(U(x,t))=Df(u(x-ct,t)).\frac{\partial}{\partial x}u(x-ct,t).$$

Therefore from (18)

(20) $$\frac{\partial}{\partial t}U(x,t)+\frac{\partial}{\partial x}(f(U(x,t))+cU(x,t))=0.$$

If we know that all velocities of waves for (18) are contained in a finite interval [*a*,*b*], then all velocities of waves for (20) are contained in [*a*+*c*,*b*+*c*], therefore they are positive if *c* is chosen larger than -*a*. If one can prove convergence of a numerical scheme for (20) when all velocities are positive, then using (19) in the reverse order would prove convergence of a numerical scheme for (18). From the maximum principle in u proved in theorem 2, this method can be applied to the system of pressureless fluid dynamics with $c \geq -\inf(u_0(x))$. After this change of variable, the new equations obtained when pressureless fluid dynamics is considered in (18) are:

(21) $$R_t+(RU)_x+cR_x=0,$$
$$(RU)_t+(RU^2)_x+c(RU)_x=0.$$

System (21) is treated with the above scheme in which the velocities are $U_i^n+c$ instead of $u_i^n$ in the auxiliary function *L*, i.e.

$$R_i^{n+1} = R_{i-1}^n\ L(-1+rU_{i-1}^n+rc,\ rU_{i-1}^n+rc) + R_i^n\ L(rU_i^n+rc, 1+rU_i^n+rc)$$
$$+R_{i+1}^n L(1+rU_{i+1}^n+rc, 2+\ rU_{i+1}^n+rc),$$

$$(RU)_i^{n+1} = (RU)_{i-1}^n\ L(-1+rU_{i-1}^n+rc,\ rU_{i-1}^n+rc) + (RU)_i^n L(rU_i^n+rc, 1+rU_i^n+rc)$$
$$+(RU)_{i+1}^n L(1+rU_{i+1}^n+rc, 2+rU_{i+1}^n+rc),$$

completed by

$$U_i^{n+1} = \frac{(RU)_i^{n+1}}{R_i^{n+1}}.$$

The CFL condition is

$$r\ max\{|U_i^n|+c\} \leq 1.$$

Instead of performing the translation (of $Nrc$ to the left if $N$ is the total number of iterations) at the end of calculations, one can choose $r$ and $c$ such that $r$ max$\{|U_i^n|\}$ $\leq 1/2$ and $rc=1/2$, and change $i$ into $i$-$1$ every two time steps. The final result from the two strategies is the same.

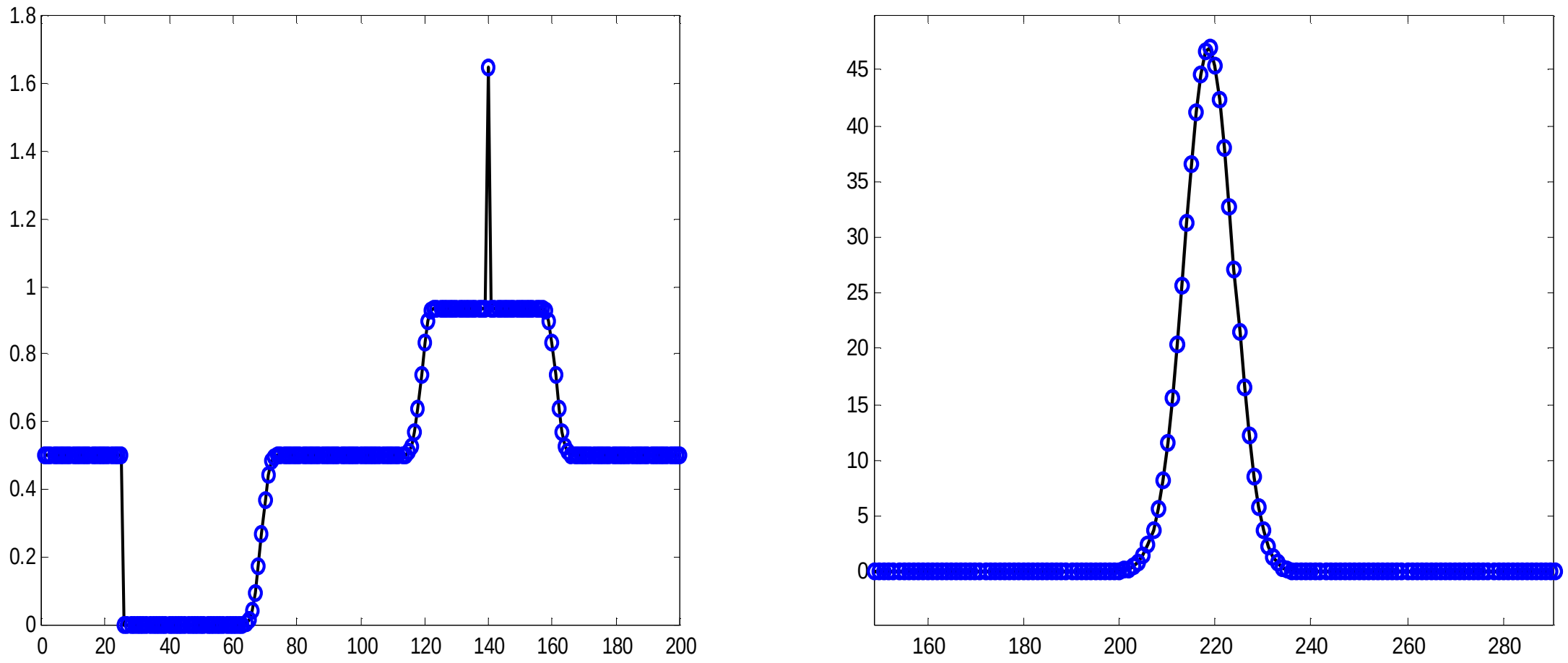

**figure 3: the previous scheme and the slightly modified one of this section give exactly same results as expected. Evidence from the tests top-left in figure 1 and bottom-right in figure 2; black 'line' = scheme of section 6, blue 'o' = scheme of this section.**

The scheme in section 6, which corresponds to $c=0$ in the formulas above, permits therefore a smaller value of r. Therefore it is better suited in numerical works and (motivated by Theorem 4) one conjectures it converges (since it is very close to the one in this section).

**Theorem 4** (convergence of the scheme to a solution of the Cauchy problem): For any initial data $\rho$ $L^1$ and $u$ $L^\infty$ the slightly modified scheme is convergent to a pair ($\rho$,$u$) of a measure $\rho$ and a $L^\infty$ function $u$ which is solution of the system of pressureless fluid dynamics in the sense of distributions on $IR^n \times [0,+\infty)$ .The result holds in any space dimension.

Proof: The maximum principle in $u$ is proved as in theorem 2. Then the convergence proof is a direct adaptation of the proof of theorem 3 in which one replaces $(\rho u)\psi_x$ by $(\rho u + c\rho)\psi_x$ (first equation) and $(\rho u^2)\psi_x$ by $(\rho u^2 + c\rho u)\psi_x$ (second equation), as well as in the corresponding discretizations of these formulas. It is clear the proof adapts to any space dimension by choosing suitable constants $c_x$, $c_y$,... in each direction.

## 10. Application. A toy-model of evolution of small perturbations in a static background fluid.

Let us consider the system of pressureless fluid dynamics coupled to gravitation ( [10] p233, [5] p207), in 1D:

$$\rho_t + (\rho u)_x = 0 ,$$

$$(\rho u)_t + (\rho u^2)_x + \rho \Phi_x = 0 ,$$

$$\Delta\Phi = 4\pi G\rho , \Phi = 0 \quad on \quad the \quad boundary .$$

This system is splitted into pressureless fluid dynamics and the system

$$\rho = \text{constant in time},$$
$$(\rho u)_t + \rho\Phi_x = 0 ,$$
$$\Delta\Phi = 4\pi G\rho ,$$

for which a numerical solution is straightforward : one computes $\Phi$, then $\Phi_x$ and integrates in t. Extension to 2D is straightforward.

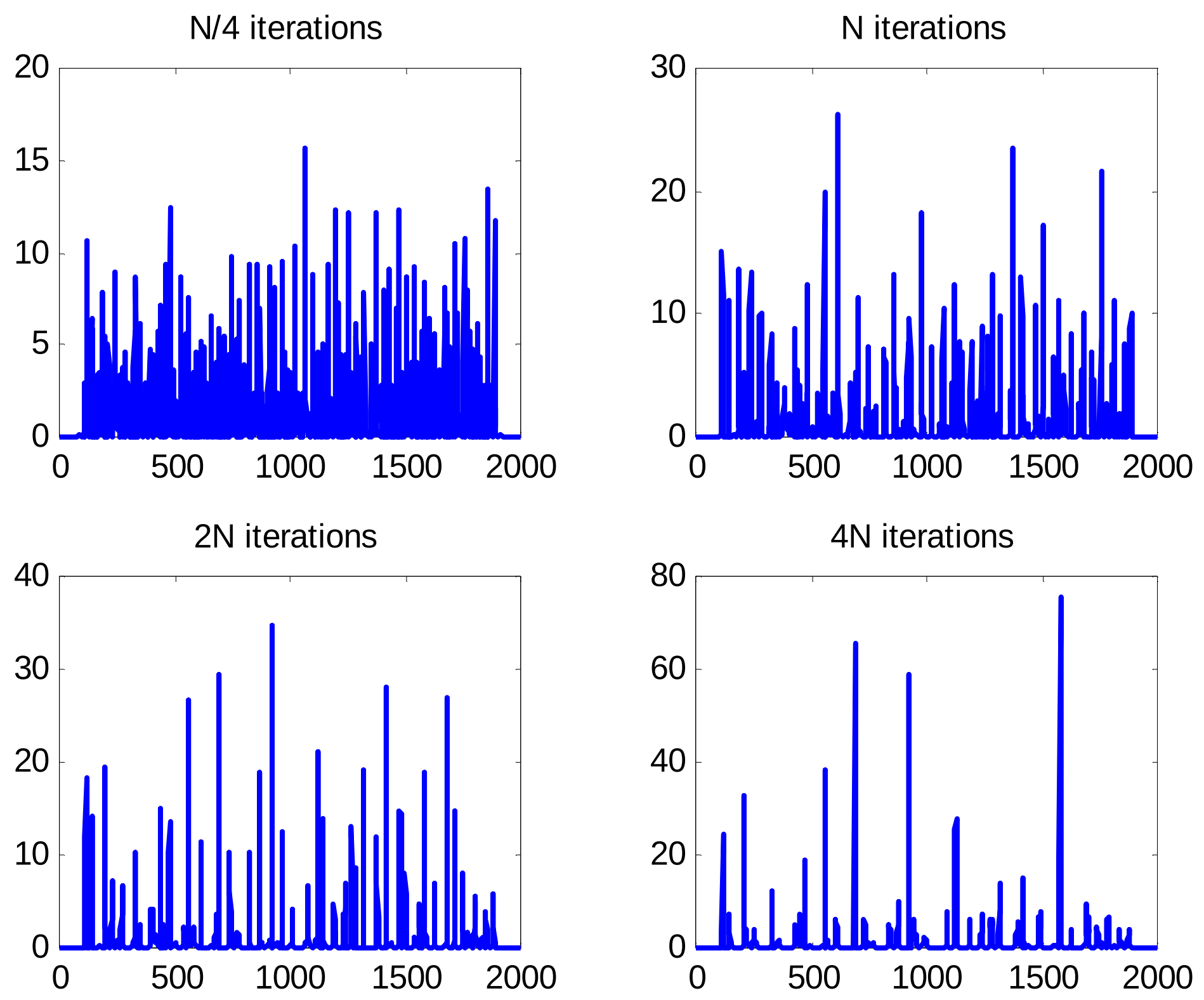

**figure 4: 1D case including gravitation. Depiction of values of the density $\rho$ at 4 successive times. The initial values of $\rho$ are chosen at random between 0.9 and 1.1. The initial values of u are chosen at random between -0.5 and 0.5. We observe the effect of gravitation. After a few time steps the matter is condensed in form of bound structures. After more time steps peaks agglomerate more and more, and grow since the total amount of matter is constant, forming larger structures with time.**

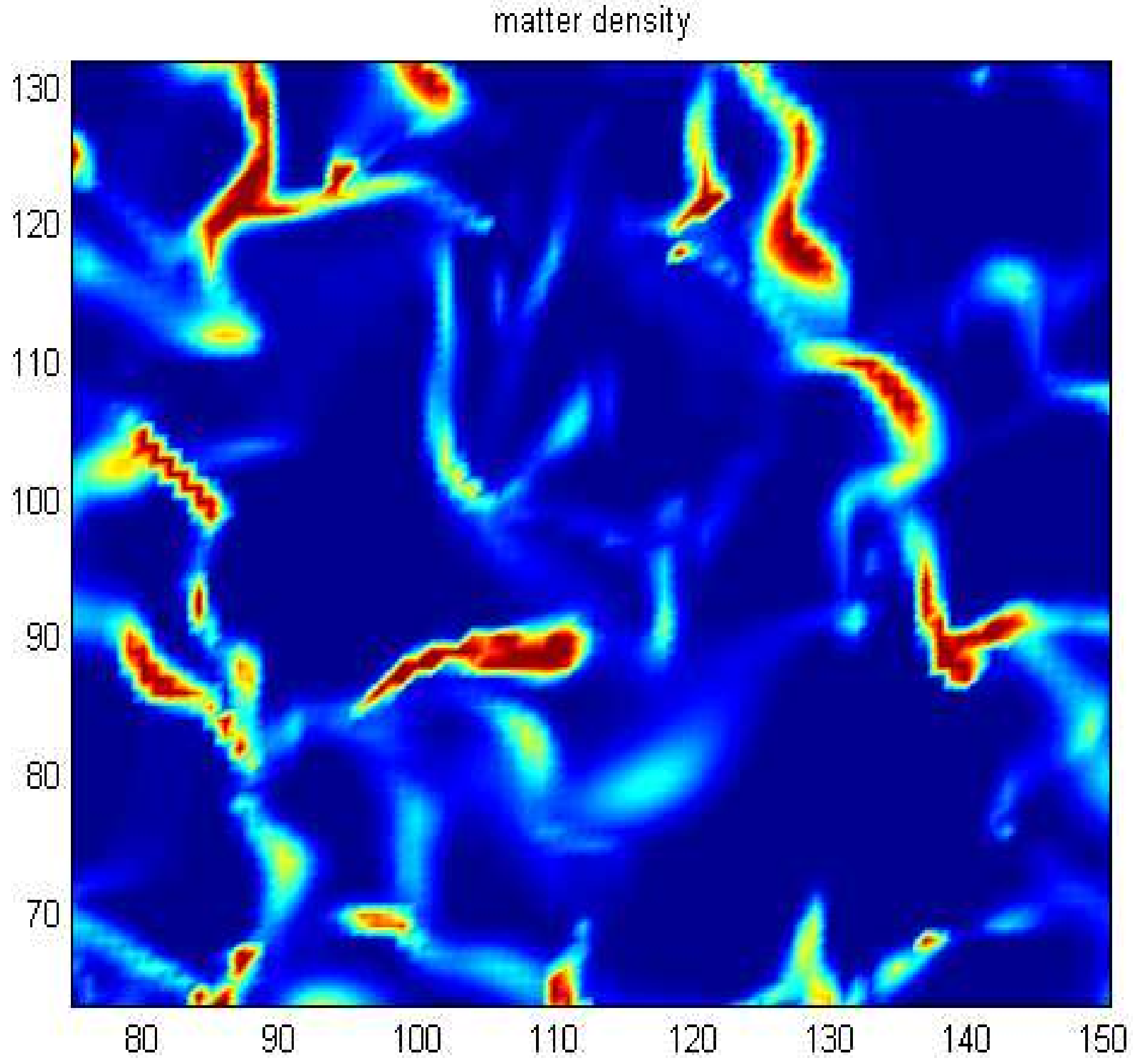

**Figure 5: 2D case including gravitation in a static background. Top: the dark-blue regions are void, the red ones are full of matter. In a vast region, one observes a patchwork of voids, clusters and filaments of matter. One observes already (in 2D and static universe) the typical structure filament-cluster-void network observed in [5] p308, p333, [9] cover, [10] cover, p490, p458. The time evolution one observes is exactly similar to the one depicted in [5] p308, which can be guessed from the 1D evolution depicted in figure 4 and the above figure. This 2D test, as well as the much more evolved ones in [6] (expanding background, relativistic, multifluid), can be done on any standard PC (here meshgrid 200$\times$200, 100 iterations, computation time 5 minutes).**

## Conclusion

In the case of numerous changes in signs of velocities the method presented in this paper permits to extend the applicability of the scheme in [1], which introduced the pioneering idea to project delta waves in a Godunov scheme. It has an interpretation from the physical viewpoint (theorem 1), mathematically strengthened by the result of stability (theorem 2) and the convergence theorems (theorems 3, 4).

The scheme has given the correct results in all numerical tests when an exact solution exists (as it could be presumed from the convergence theorem 4). In particular, the scheme in this paper works well for problems where the velocity changes sign in regions where the density varies smoothly (figure 3). It shows a good accuracy, taking into account it is very fast and with a rather large CFL condition, but it is not free of numerical dissipation (to be possibly diminished by numerical techniques, but this appears useless in applications to cosmology [6] since it has been observed gravitation can reduce the support of a delta peak to one cell). It extends at once (without dimensional splitting) into 2D and 3D schemes, whose use is easy inside the technique of splitting of equations [1,6] for numerical approximation of more complicated systems that cannot be conveniently treated globally.

Systems more complicated than pressureless fluid dynamics can be treated in view of a Godunov scheme by the original method of projection of delta peaks presented in this paper. As an application, the Newtonian and relativistic systems that describe matter and energy submitted to gravitation in an expanding background are solved numerically in 1D, 2D and 3D in [6], as well as multi-fluid extensions. The numerical tests in [6] provide a 2D-description of some of the main events of the theory of large structure formation in cosmology, that explains the structure of galaxies and galaxy clusters observed in the today universe.

## Appendix 3: Calculations of the delta peak solutions of the Riemann problem for the system of pressureless fluid dynamics

The results can also be found in[3,7]. We seek a solution of the form

$$u(x,t) = u_l + \Delta u.\ H(x-ct)\ ,$$
$$\rho(x,t) = \rho_l + \Delta\rho.H(x-ct) + \alpha t\delta(x-ct),$$

where $H$=Heaviside function, $\delta$=Dirac delta function. We insert these formulas into the system

$$\rho_t + (\rho u)_x = 0\ ,$$
$$(\rho u)_t + (\rho u^2)_x = 0\ .$$

After separation of the terms in factor of $t$ and without $t$ the first equation gives the two equations

(22) $$-c\Delta\rho + \alpha + \Delta(\rho u) = 0,\ c\delta = u_l\delta + \Delta u H\delta\ ,$$

Taking into account the line above, the term without t from the second equation appears as

$$-c\Delta(\rho u) + c\alpha + \Delta(\rho u^2) = 0 .$$

This equation and the first equation in (22) give

$$(\Delta\rho)c^2 - 2\Delta(\rho u)c + \Delta(\rho u^2) = 0.$$

Solving this equation in c

$$c = \frac{\Delta(\rho u) \pm \Delta u\sqrt{\rho_r \rho_l}}{\Delta\rho},$$

which, in turn, gives

$$\alpha = \pm\Delta u\sqrt{\rho_l \rho_r} .$$

If $\Delta u < 0$, i.e. $u_l > u_r$, one has a peak in density, therefore $\alpha > 0$ and thus

$$\alpha = -\Delta u\sqrt{\rho_l \rho_r}, c = \frac{\Delta(\rho u) - \sqrt{\rho_l \rho_r}}{\Delta\rho} = \frac{u_l\sqrt{\rho_l} + u_r\sqrt{\rho_r}}{\sqrt{\rho_l} + \sqrt{\rho_r}} .$$

If $\Delta u > 0$, i.e. $u_l < u_r$, the above formulas give a negative density peak: indeed one has a void region according to (5). The negative density peak depicted by the above formula is a nonphysical solution. One can check the other choices of signs in $\pm$ are nonsense. The formula for $\beta$ in (3) is obtained by multiplying $\rho$ and *u:*

$$\beta\delta = u_l\alpha\delta + (\Delta u)\alpha H\delta .$$

Integration and use of (22) give:

$$\beta = u_l\alpha + \alpha(c - u_l) = \alpha c .$$